\documentclass[a4paper]{amsart}
\usepackage{amsmath, amssymb, amsfonts}
\numberwithin{equation}{section}
\theoremstyle{plain}
 \newtheorem{theorem}{Theorem}[section]
 \newtheorem{corollary}[theorem]{Corollary}
 \newtheorem{lemma}[theorem]{Lemma}
 \newtheorem{proposition}[theorem]{Proposition}
\theoremstyle{definition}
 \newtheorem{definition}[theorem]{Definition}
 \newtheorem{example}[theorem]{Example}
\theoremstyle{remark}
 \newtheorem{remark}[theorem]{Remark}
\def\ord{\operatorname{ord}}
\def\Sol{\operatorname{Sol}}
\DeclareMathOperator{\sh}{sinh}

\def\Ad{\operatorname{Ad}}
\def\sign{\operatorname{sgn}}
\def\Hom{\operatorname{Hom}}

\def\p{\partial}
\begin{document}
\title[Commuting differential operators
with regular singularities]{Commuting differential operators\\
with regular singularities}%
\author{Toshio Oshima}
\address{Graduate School of Mathematical Sciences,
University of Tokyo, 7-3-1, Komaba, Meguro-ku, Tokyo 153-8914, Japan}
\email{oshima@ms.u-tokyo.ac.jp}
\maketitle
\begin{abstract}
We study a system of partial differential equations defined by
commuting family of differential operators with regular singularities.
We construct ideally analytic solutions depending on a holomorphic parameter.
We give some explicit examples of differential operators related 
to $SL(n,\mathbb R)$ and completely integrable quantum systems.
\end{abstract}
\maketitle
\section{Introduction}
The invariant differential operators on a semisimple symmetric space
have regular singularities along the boundaries of the space which is realized 
in a compact manifold by \cite{OR2}.
In the case of a Riemannian symmetric space $G/K$,  the study of such operators
in \cite{KO} enables \cite{K--} to have the Poisson integral expression
of any simultaneous eigenfunction of the operators.
Here $G$ is a connected real semisimple Lie group with finite center
and $K$ is its maximal compact subgroup.

In the case of the group manifold $G$, which is an example of a semisimple 
symmetric space, Harish-Chandra gives an asymptotic expansion of a right 
and left $K$-finite eigenfunction, which plays an important role in the 
harmonic analysis on $G$ (cf.~\cite{Ha}).
He uses only the Casimir operator to get the asymptotic expansion, which
suggests us that one operator controls other operators together with some 
geometry.

On the other hand, the Schr\"odinger operator corresponding to 
Calogero-Moser-Sutherland system with a trigonometric potential function 
(cf.~\cite{Su})
or a Toda finite chain (cf.~\cite{To}) is completely 
integrable and the integrals with higher orders
are uniquely characterized by the Schr\"odinger operator and so are the 
simultaneous eigenfunctions.
These integrals also have regular singularities at infinity.

In this note we study a general commuting system of differential operators
with regular singularities by paying attention to the fact that 
an operator characterizes the system.
Our argument used in this note is based on expansions in power series and
hence it is rather elementary compared to that in \cite{KO} and \cite{OV2} 
where a microlocal method is used.

In fact we will study matrices of differential operators which may not
commute with others in the system but satisfy a certain condition because 
it is better to do so even in the study of commuting scalar differential 
operators.
Some of its reasons will be revealed in the proof of Theorem~\ref{thm:induce},
that of Theorem~\ref{thm:nopara}, Remark~\ref{rem:nonsplit} ii) etc.

In \S\ref{sec:CD} we study differential operators which commute one
operator. 
We will see that the \emph{symbol map} $\sigma_*$ plays an important role. 
In the case of the first example above the map corresponds to Harish-Chandra's
isomorphism of the invariant differential operators.
In the case of the Schr\"odinger operator above it corresponds to
the commutativity among the integrals with higher orders.

In \S\ref{sec:SI} we construct some of multivalued holomorphic solutions 
of the system around the singular points which we call \emph{ideally analytic 
solutions} and then
in \S\ref{sec:Ind} we study the \emph{induced equations} of other 
operators, which assures that the solutions automatically satisfy some other 
differential equations.

In \S\ref{sec:HS} we study the holonomic system of differential equations
with constant coefficients holomorphically depending on a parameter, which 
controls the \emph{leading terms} of the ideally analytic solutions.

In \S\ref{sec:AI} we study a \emph{complete system of differential equations
with regular singularities}
which means that the system is sufficient to formulate a boundary value problem
along the singularities and we describe all the ideally analytic solutions.
In particular, when the system has a holomorphic parameter, we
construct solutions depending holomorphically on the parameter.
It is in fact useful to introduce a parameter for the study of a specific 
system by holomorphically deforming it to generic simpler ones.

In \S\ref{sec:SL} and \S\ref{sec:CI}\ we give some explicit examples of 
the systems related to $SL(n,\mathbb R)$ and the completely
integrable quantum systems with regular singularities at infinity,
respectively.
Moreover we give Theorem~\ref{thm:CI2} in the case of completely integrable
quantum systems with two variables.

\section{Commuting differential operators with regular singularities}
\label{sec:CD}
For a positive integer $m$ and a ring $R$
we will denote by $M(m,R)$ the ring of square matrices 
of size $m$ with components in $R$ and by $R[\xi]$ the ring of polynomials 
of $n$ indeterminates $\{\xi_1,\ldots,\xi_n\}$ if $\xi=(\xi_1,\ldots,\xi_n)$.
The $(i,j)$-component of $A\in M(m,R)$ is denoted by $A_{ij}$ and we 
naturally identify $M(1,R)$ with $R$.

Let $M$ be an $(n+n')$-dimensional real analytic manifold and let $N_i$ be 
one-codimensional submanifolds of $M$ such that $N_1,\ldots,N_n$ are normally
crossing at $N=N_1\cap\cdots\cap N_n$.
We assume that $M$ and $N$ are connected.
We will fix a local coordinate system 
$(t,x)=(t_1,\ldots,t_n,x_1,\ldots,x_{n'})$ around a point $x^o\in N$ 
so that $N_i$ are defined by the equations $t_i=0$, 
respectively.  

Let $\mathcal A_N$ denote the space of real analytic functions on $N$
and $\mathcal A_M$ the space of real analytic functions defined on a 
neighborhood of $N$ in $M$.
For $\alpha=(\alpha_1,\ldots,\alpha_n)$, 
$\beta=(\beta_1,\ldots,\beta_n)\in\mathbb Z^n$ we put
\begin{align*}
 |\alpha| &= \alpha_1+\cdots+\alpha_n,\\
 \alpha<\beta\ &\Leftrightarrow\ \alpha_i\le\beta_i
 \quad\text{for }i=1,\ldots,n
 \text{ and }\alpha\ne\beta.
\end{align*}

Let $\mathbb N$ be the set of non-negative integers.
We will denote
\begin{equation*}
\begin{cases}
 \vartheta_i = {t_i}\frac{\p}{\p t_i},\quad
 \p_x = \left(\frac{\p}{\p x_1},\ldots,\frac{\p}{\p x_{n'}}\right),\\
 \vartheta^\alpha = \vartheta_1^{\alpha_1}\cdots\vartheta_n^{\alpha_n}
 \quad\text{for } \alpha=(\alpha_1,\ldots,\alpha_n)\in\mathbb N^n,\\
 \p_x^\beta = \frac{\p^{|\beta|}}{\p x_1^{\beta_1}\cdots 
\p x_{n'}^{\beta_{n'}}}
 \quad\text{for } \beta=(\beta_1,\ldots,\beta_{n'})\in\mathbb N^{n'},\\
 t^\lambda = t_1^{\lambda_1}\cdots t_n^{\lambda_n}
 \quad\text{for } \lambda=(\lambda_1,\ldots,\lambda_n)\in\mathbb C^n.
\end{cases}
\end{equation*}

Let $\mathcal D_M$ and $\mathcal D_N$ denote the rings of differential 
operators on $M$ and $N$ with coefficients in $\mathcal A_M$ and
$\mathcal A_N$, respectively.
\begin{definition}
Let $\widetilde{\mathcal D}_*$ denote the subring of $\mathcal D_M$ whose
elements $P$ have the form
\begin{equation}\label{eq:P0}
  P = \sum_{(\alpha,\beta)\in\mathbb N^{n+n'}}
        a_{\alpha,\beta}(t,x)
      \vartheta^\alpha\p_x^\beta\text{ \ with \ }
  a_{\alpha,\beta}(t,x)\in\mathcal A_M.
\end{equation}
Here the sum above is finite.
Moreover $\mathcal D_*$ denotes the subring of $\widetilde{\mathcal D}_*$
whose elements $P$ of the form \eqref{eq:P0} satisfy
\begin{equation}\label{eq:P1}
   a_{\alpha,\beta}(0,x)=0\quad \text{if \ }\beta\ne0.
\end{equation}
When $P$ is an element of $\mathcal D_*$, $P$ is said to have 
\emph{regular singularities in the weak sense along the set of walls 
$\{N_1,\ldots,N_n\}$ with the edge $N$ (cf.~\cite{KO})}.

Let define a map $\sigma_*$ of $\widetilde{\mathcal D}_*$ to $\mathcal D_N[\xi]$ by
\begin{equation*}
 \sigma_*(P)(x,\xi,\p_x) :=\sum_{\alpha\in\mathbb N^n,\,\beta\in\mathbb N^{n'}} 
 a_{\alpha,\beta}(0,x)\xi^\alpha \p_x^\beta
\end{equation*}
for $P$ in \eqref{eq:P0}. 
Then
\begin{equation*}
 t^{-\lambda}Pt^{\lambda}\phi(t,x)\bigm|_{t=0}=\sigma_*(P)(x,\lambda,\p_x)
 \phi(0,x)
 \quad\text{for }\phi\in\mathcal A_M\text{ and }\lambda\in\mathbb C^n.
\end{equation*}
\end{definition}

Here we note that the condition $P\in\widetilde{\mathcal D}_*$ equals
\[
  t^{-\lambda}Pt^{\lambda}\phi(t,x)\in \mathcal A_M
  \text{ \ for \ }\forall\phi(t,x)\in\mathcal A_M
\]
and $\sigma_*$ is a ring homomorphism of $\widetilde{\mathcal D}_*$ to 
$\mathcal D_N[\xi]$ and $\sigma_*(\mathcal D_*)=\mathcal A_N[\xi]$.

For $k\in\mathbb N$ and $P\in\widetilde{\mathcal D}_*$ with the form 
\eqref{eq:P0} we put
\begin{equation*}
     \sigma_k(P)(t,x,\xi,\tau) := \sum_{|\alpha|+|\beta|=k}
                    a_{\alpha,\beta}(t,x)\xi^\alpha\tau^\beta
\end{equation*}
and then the order of $P$, which is denoted by $\ord P$, is the maximal 
integer $k$ with $\sigma_k(P)\ne0$.

For $P=\Bigl(P_{ij}\Bigr)_{\substack{1\le i\le m\\ 1\le j\le m}}
\in M(m,\widetilde{\mathcal D}_*)$, the order of $P$ is 
defined to be the maximal order of the components of $P$ and denoted by 
$\ord P$.  We put
\begin{align*}
  \sigma(P)&:=
   \Bigl(\sigma_{\ord P}(P_{ij})\Bigr)_{\substack{1\le i\le m\\ 1\le j\le m}}
   \in M(m,\mathcal A_M[\xi,\tau]),\\
  \sigma_*(P)&:= 
   \Bigl(\sigma_*(P_{ij})\Bigr)_{\substack{1\le i\le m\\ 1\le j\le m}}
   \in M(m,\mathcal D_N[\xi]),\\
  \bar\sigma_*(P)&:=
  \sigma(P)(0,x,\xi,\p_x)\in  M(m,\mathcal D_N[\xi]).
\end{align*}
Then as a polynomial of $\xi$, 
$\bar\sigma_*(P)$ is the homogeneous part of $\sigma_*(P)$ whose degree
equals $\ord P$.
For $P,\,Q\in\widetilde{\mathcal D}_*$, we note that
$\sigma(PQ)=\sigma(P)\sigma(Q)$ and 
\begin{equation*}
 \begin{split}
 \sigma_{\ord P+\ord Q-1}([P,Q])
 &=\sum_{i=1}^n
    \Bigl(\frac{\p \sigma(P)}{\p \xi_i}t_i\frac{\p \sigma(Q)}{\p t_i}
   -\frac{\p \sigma(Q)}{\p \xi_i}t_i\frac{\p \sigma(P)}{\p t_i}\Bigr)\\
 &\ +\sum_{j=1}^{n'}
    \Bigl(\frac{\p \sigma(P)}{\p \tau_j}\frac{\p\sigma(Q)}{\p x_j}
   -\frac{\p \sigma(Q)}{\p \tau_j}\frac{\p\sigma(P)}{\p x_j}\Bigr).
 \end{split}
\end{equation*}
\begin{theorem}\label{thm:com}
Let $P$ and $Q$ be nonzero elements of $M(m,\widetilde{\mathcal D}_*)$
such that $[P,Q]=0$, $P\in M(m,\mathcal D_*)$ and $\sigma(P)$ is a scalar 
matrix satisfying
\begin{equation}
 \sum_{\nu=1}^n\gamma_\nu\frac{\p\bar\sigma_*(P)}{\p\xi_\nu}\not\equiv 0
 \quad\text{for \ }\forall\gamma\in\mathbb N^n\setminus\{0\}.\label{eq:nondeg}
\end{equation}
Here ``$\not\equiv0$" means ``not identically zero".
Suppose that $\sigma_{\ord P-1}(P)$ or $\sigma(Q)$ is a scalar matrix.
Then $[\sigma_*(P),\sigma_*(Q)]=0$ and $\bar\sigma_*(Q)\ne0$.
Moreover if $\sigma(P)(t,x,\xi,\tau)$ does not depend on $t$, 
so does $\sigma(Q)(t,x,\xi,\tau)$.
\end{theorem}
\begin{proof}
Since $\sigma_*$ is an algebra homomorphism,
$[\sigma_*(P),\sigma_*(Q)]=\sigma_*([P,Q])=0$.

Put $r_P=\ord P$ and $r_Q=\ord Q$.
Fix $i$ and $j$ such that $\sigma_{r_Q}(Q_{ij})\ne 0$.
Note that the assumption implies
\[ 
  \sigma_{r_P+r_Q-1}([P,Q]_{ij}) = \sigma_{r_P+r_Q-1}([P_{11},Q_{ij}]).
\]
Put
\begin{gather*}
 \sigma_{r_P}(P_{11})  = 
  \sum_{\substack{\beta,\,\gamma\\ |\beta|\le r_P}}
  p_{\beta,\gamma}(x,\xi)t^\gamma\tau^\beta,\ 
 \sigma_{r_Q}(Q_{ij})  = 
  \sum_{\substack{\beta,\,\gamma\\ |\beta|\le r_Q}}
  q_{\beta,\gamma}(x,\xi)t^\gamma\tau^\beta,\allowdisplaybreaks\\
 \sigma_{r_P+r_Q-1} ([P_{11},Q_{ij}]) = 
  \sum_{\substack{\beta,\,\gamma\\ |\beta|\le r_P+r_Q-1}}
  s_{\beta,\gamma}(x,\xi)t^\gamma\tau^\beta
\end{gather*}
and choose $(\beta^o,\gamma^o)\in\mathbb N^{n'+n}$ such that
\begin{equation*}
 \begin{cases}
   q_{\beta^o,\gamma^o}\ne 0,\\
   q_{\beta,\gamma}=0&\text{if } \gamma < \gamma^o,\\
   q_{\beta,\gamma^o}=0&\text{if } \beta > \beta^o.
 \end{cases}
\end{equation*}
Then
\begin{equation}\label{eq:uniq}
 s_{\beta^o,\gamma^o}t^{\gamma^o}\tau^{\beta^o}=\Bigl(\sum_{\nu=1}^n
  \frac{\p p_{0,0}}{\p\xi_\nu}\gamma^o_\nu\Bigr)
  \Bigl(q_{\beta^o,\gamma^o}t^{\gamma^o}\tau^{\beta^o}\Bigr),
\end{equation}
which proves
the first claim in the theorem because the condition $[P,Q]=0$ 
with the assumption of the theorem means $\gamma^o=0$.

Moreover suppose $p_{\beta,\gamma}=0$ for $\gamma\ne0$.  
Then \eqref{eq:uniq} is valid for any $\gamma^o\in\mathbb N^n$
and $\beta^o\in\mathbb N^{n'}$ satisfying $q_{\beta,\gamma^o}=0$ for 
$\beta>\beta^o$
and hence the condition $[P,Q]=0$ means
$q_{\beta^o,\gamma^o}=0$ if $\gamma^o\ne0$.
Thus $q_{\beta,\gamma^o}=0$ if $\gamma^o\ne0$.
\end{proof}
\begin{corollary}\label{cor:commute}
Let $P\in M(m,\mathcal D_*)$
such that $\sigma(P)$ and  $\sigma_{\ord P-1}(P)$ are scalar matrices.
Suppose $\bar\sigma_*(P)$ satisfies \eqref{eq:nondeg}.
Then the map
\begin{equation*}
 \begin{split}
  \sigma_*:\,M(m,\widetilde{\mathcal D}_*)^P:= \{Q\in M(m,\widetilde{\mathcal D}_*)\,;
 \,[P,Q]=0\}
  &\to M(m,\mathcal D_N[\xi]),\\
 Q&\mapsto \sigma_*(Q)
 \end{split}
\end{equation*}
is an injective algebra homomorphism.

In particular, when $m=1$, 
$\mathcal D_*^P:=\{Q\in\mathcal D_*\,;\,[P,Q]=0\}$ is commutative.
\end{corollary}
\begin{proof}\label{cor:HC}
Since $\sigma_*$ is an algebra homomorphism and the condition
$Q_1,\,Q_2\in M(m,\widetilde{\mathcal D}_*)^P$ implies
$[Q_1,Q_2]\in M(m,\widetilde{\mathcal D}_*)^P$, this corollary is a direct consequence 
of Theorem~\ref{thm:com}.
\end{proof}
\begin{remark} {\rm i)}
Retain the notation in Theorem~\ref{thm:com}.
Then \eqref{eq:nondeg} is valid for $P\in M(m,\mathcal D_*)$ if
$n$ functions
$\frac{\p \bar\sigma_*(P)}{\p\xi_1},\ldots,
\frac{\p \bar\sigma_*(P)}{\p\xi_n}$
are linearly independent over $\mathbb R$.
In particular, if $\ord P=2$ and $\bar\sigma_*(P)$ is a scalar matrix,
the condition that
\begin{equation*}
 \text{the matrix }
 \left(\frac{\p^2\bar p}{\p\xi_i\p\xi_j}\right)
 _{\substack{1\le i\le n\\1\le j\le n}}
 \quad\text{is invertible for generic }x\in N
\end{equation*}
implies \eqref{eq:nondeg}.
Here $\bar p$ is the diagonal element of $\bar\sigma_*(P)$.

{\rm ii)}
The assumption $P\in M(m,\mathcal D_*)$ is necessary in Theorem~\ref{thm:com}.
For example, 
$[t\frac\p{\p t}+x\frac \p{\p x}, t\frac\p{\p x}]=0$ and $\sigma_*(t\frac\p{\p x})=0$.
Moreover we note that
\[
\left[
\begin{pmatrix}
 t\frac\p{\p t} &  0\\
 0              & t\frac\p{\p t}+1
\end{pmatrix},
\begin{pmatrix}
 0 & t\\
 0 & 0
\end{pmatrix}
\right]=0.
\]
This gives an example such that $\sigma_{\ord P-1}(P)$ and
$\sigma(Q)$ are not scalar matrices.

{\rm iii)} The invariant differential operators on a Riemannian symmetric space
$G/K$ of non-compact type have regular singularities along the boundaries of
a realization of the space constructed by \cite{OR1} and the map $\sigma_*$
of $\mathcal D_*^P$ to $\mathcal A_N[\xi]$ in Corollary~\ref{cor:HC}
corresponds to Harish-Chandra isomorphism (cf.~\cite{K--}).

The element of the universal enveloping algebra $U(\mathfrak g)$ of 
the Lie algebra of $G$ defines a differential operator on the realization 
of $G/K$ through the infinitesimal action of the left translation by elements
of $G$.
Then the differential operator is an element of $\widetilde{\mathcal D}_*$.

Moreover the invariant differential operators on a semisimple symmetric
space whose rank is larger than its real rank are in $\widetilde{\mathcal D}_*$
(cf.~\cite{OV2}).

The radial parts of the Casimir operator acting on $K$-finite sections of
certain homogeneous vector bundle of $G$ satisfy the assumption of
Theorem~\ref{thm:com} (cf.~\eqref{eq:Sph} and \eqref{eq:KWhit} for examples).
\end{remark}
\section{Ideally analytic solutions without logarithmic terms}\label{sec:SI}
For a subset $\Sigma$ of $\mathbb N^n$ define
\begin{align*}
 \overline\Sigma:&=\bigl\{\alpha\in\mathbb N^n\,;\,
 \{\alpha+\gamma\,;\,\gamma\in \mathbb N^n\}\cap\Sigma\ne\emptyset\bigr\},\\
 \p\Sigma:&=\bigl\{\alpha\in\mathbb N^n\setminus\overline\Sigma\,;\,
   \text{ there exists }\gamma\in\bar\Sigma\text{ such that }
   \sum_{i=1}^n|\alpha_i-\gamma_i|=1\}.
\end{align*}
Moreover we denote by $\hat{\mathcal A}_M$ the ring of formal power series
of $t=(t_1,\ldots,t_n)$ with coefficients in $\mathcal A_N$.

\begin{theorem}\label{thm:sol}

Let $P\in M(m,\mathcal D_*)$.

{\rm i)}
Let $\Sigma$ be a subset of $\mathbb N^n$ such that
\begin{equation*}
  \det\bigl(\sigma_*(P)(x,\gamma)\bigr)\not\equiv 0
  \quad\text{for \ }\forall\gamma\in\mathbb N^n\setminus\Sigma.
\end{equation*}
Let 
$\hat u(t,x)=\sum_{\alpha\in\mathbb N^n}
u_{\alpha}(x)t^\alpha\in\hat{\mathcal A}_M^m$
be a formal solution of $P\hat u=0$.
Then $\hat u=0$ if $u_\alpha=0$ for $\forall\alpha\in\Sigma$.

Hereafter in this theorem suppose
\begin{align}\label{eq:converg}
 \det \bar\sigma_*(P)(x,\xi)\ne0
 \quad\text{for }\forall\xi=(\xi_1,\ldots,\xi_n)\in [0,\infty)^n\setminus\{0\}
 \ \text{and } \forall x\in N.
\end{align}

{\rm ii)} If $\hat u\in\hat{\mathcal A}_M^m$ satisfies 
$P\hat u\in\mathcal A_M^m $, then $\hat u\in\mathcal A_M^m$.

{\rm iii)}
Fix $f\in\mathcal A_M^m$, a point $x^o\in N$ and a finite subset $\Sigma$ of 
$\mathbb N^n$ such that
\begin{equation*}
  \det\bigl(\sigma_*(P)(x^o,\gamma)\bigr)\ne 0
  \quad\text{for \ }\forall\gamma\in\mathbb N^n\setminus\Sigma.
\end{equation*}
By shrinking $M\ni x^o$ if necessary and denoting
\begin{align*}
 \mathcal A_M(P^{-1}f)&:=
  \{u\in \mathcal A_M^m\,;\,Pu=f\},
\allowdisplaybreaks\\
 \mathcal A_M(P^{-1}f)^{\Sigma}&:=
  \{\bar u=\sum_{\alpha\in\overline\Sigma}u_\alpha(x)t^{\alpha}
  \in \mathcal A_M^m\,;\,
  P\bar u \equiv f \mod
 \sum_{\beta\in\p\Sigma}\mathcal A_M^m t^\beta\},
\end{align*}
the natural restriction map
\begin{equation*}
 \mathcal A_M(P^{-1}f)\, \xrightarrow{\sim}\, \mathcal A_M(P^{-1}f)^{\Sigma},\quad
 \sum_{\alpha\in\mathbb N^n}
 u_\alpha(x)t^{\alpha}\mapsto \sum_{\alpha\in\overline\Sigma}
 u_\alpha(x)t^{\alpha}
\end{equation*}
is a bijection. Here in particular
\begin{equation*}
 \mathcal A_M(P^{-1}f)^{\{0\}} 
 =\{u\in\mathcal A_N^m\,;\,
   \sigma_*(P)\bigl(x,0)u = f|_{t=0}\}.
\end{equation*}
\end{theorem}
\begin{proof}
The proof proceeds in a similar way as in \cite[Theorem~2.1]{OV1} where we 
studies the same problem with $n=1$.

We may assume $x^o=0$.
Expanding functions in convergent power series of $(t,x)$ at 
$(0,0)$, we will prove the theorem in a neighborhood of 
$(0,0)$.

Put $r=\ord P$ and
\begin{equation*}
 \begin{split}
  P = \sigma_*(P)(x,\vartheta)
    +\sum_{\substack{(\alpha,\beta)\in\mathbb N^{n+n'}\\
    |\alpha|+|\beta|\le r}}p_{\alpha,\beta}(t,x)\vartheta^\alpha\p_x^\beta.
 \end{split}
\end{equation*}
Then 
$
  p_{\alpha,\beta}(0,x)=0.
$
For a finite subset $\Sigma\subset\mathbb N^n$ and 
\begin{equation*}
  \hat u(t,x) = \sum_{\alpha\in\mathbb N^n}
  \hat u_\alpha(x)t^\alpha\in\hat{\mathcal A}_M^m,
\end{equation*}
put
\begin{equation*}
 \bar u(t,x) = \sum_{\alpha\in\overline\Sigma}\hat u_\alpha(x)t^\alpha.
\end{equation*}

Suppose 
$P\hat u\equiv f\mod\sum_{\alpha\in\p\Sigma}\hat{\mathcal A}_M^m t^\alpha$.
Put $h=f-P\bar u$.
Then
\begin{equation*}
 h=\sum_{\alpha\in\mathbb N^n\setminus\overline\Sigma}h_\alpha(x)t^\alpha
  = \sum_{\alpha\in\mathbb N^n\setminus\overline\Sigma,
    \ \beta\in\mathbb N^{n'}}
  h_{\alpha,\beta}t^\alpha x^\beta
 \in \mathcal A_M^m
\end{equation*}
and
\begin{equation*}
 P\hat u = f\ \Leftrightarrow\ Pu=h\quad\text{with }u=\hat u-\bar u.
\end{equation*}

Then the equation $P\hat u=f$ is equal to
\begin{equation*}
\begin{split}
  \sigma_*(P)(x,\vartheta)u &= h
  -\sum_{\substack{(\alpha,\beta)\in\mathbb N^{n+n'}\\
     |\alpha|+|\beta|\le r}}
    p_{\alpha,\beta}(t,x)\vartheta^\alpha\p_x^\beta u.\\
  u &= \sum_{\alpha\in\mathbb N^n}u_\alpha(x)t^\alpha
  \quad\text{with }
    u_\alpha(x) =
    \begin{cases}
      0&\text{for }\alpha\in\overline\Sigma,\\
      \hat u_\alpha(x)&\text{for }\alpha\in\mathbb N^n\setminus\overline\Sigma,
    \end{cases}
\end{split}
\end{equation*}
which also equals
\begin{equation}\label{eq:zen0}
\begin{split}
&\sigma_*(P)(x,\alpha^o)u_{\alpha^o}(x) = 
 h_{\alpha^o}(x) \\ 
&\quad - \text{Coef($t^{\alpha^o}$) of }
   \Bigl(\sum_{\substack{(\alpha,\beta)\in\mathbb N^{n+n'}\\ |\alpha|+|\beta|\le r}}
   p_{\alpha,\beta}(t,x)\vartheta^\alpha\p_x^\beta\Bigr)\Bigl(
     \sum_{\substack{\alpha\in\mathbb N^n\\ |\alpha|<|\alpha^o|}}
     u_{\alpha}(x)t^\alpha\Bigr)
\end{split}
\end{equation}
for $\forall\alpha^o\in\mathbb N^n\setminus\overline\Sigma$.
Here ``Coef($t^{\alpha^o}$)" means ``the coefficient of  $t^{\alpha^o}$".

Since $\det\sigma_*(P)(x,\gamma)\ne 0$ for $\gamma\in\mathbb N^n\setminus\overline\Sigma$, 
$u_{\alpha^o}(x)$ is inductively determined by $h$.

On the other hand, putting $h=0$, 
it is clear that the claim i) follows from the 
induction proving $u_{\alpha^o}=0$ by \eqref{eq:zen0} for 
$\forall\alpha^o\in\mathbb N^n\setminus\Sigma$.

Put
\begin{equation*}
 u_\alpha(x)
   = \sum_{\beta\in\mathbb N^{n'}}
     u_{\alpha,\beta}x^\beta\quad\text{with \ }u_{\alpha,\beta}\in\mathbb C.
\end{equation*}
The equation \eqref{eq:zen0} equals
\begin{align*}
 &\sigma_*(P)(0,\alpha^o)u_{\alpha^o,\beta^o}\\
 &
  = h_{\alpha_o,\beta}+ \text{Coef($x^{\beta^o}$) of }
   \bigl(\sigma_*(P)(0,\alpha^o)-\sigma_*(P)(x,\alpha^o)\bigr)
   \Bigl(\sum_{|\beta|<|\beta^o|}u_{\alpha^o,\beta}x^{\beta}\Bigr)\notag\\
 &\quad
  -  \text{Coef($t^{\alpha^o}x^{\beta^o}$) of }
   \Bigl(\sum_{\substack{(\alpha,\beta)\in\mathbb N^{n+n'}\\ 
   |\alpha|+|\beta|\le r}}
   p_{\alpha,\beta}(t,x)\vartheta^\alpha\p_x^\beta\Bigr)\Bigl(
     \sum_{\substack{(\alpha,\beta)\in\mathbb N^{n+n'}
     \\ |\alpha|<|\alpha^o|}}
     u_{\alpha,\beta}t^\alpha x^\beta\Bigr)\notag
\end{align*}
for any $\alpha^o\in\mathbb N^n\setminus\overline\Sigma$ 
and $\beta^o\in\mathbb N^{n'}$.
Hence the elements $u_{\alpha^o,\beta^o}$ of $\mathbb C^m$ satisfying this 
equation are uniquely and inductively 
determined in the lexicographic order of $(|\alpha^o|,|\beta^o|)$.
Thus to complete the proof we have only to prove 
that $\sum u_{\alpha,\beta}t^\alpha x^\beta$ is a convergent power series. 
Here we may assume $\overline\Sigma\ni\{0\}$.

In general, for formal power series 
$\psi=\sum a_{\alpha,\beta}t^\alpha x^\beta$
and $\phi=\sum b_{\alpha,\beta}t^\alpha x^\beta$ we denote
$\psi\ll\phi$ if $|a_{\alpha,\beta}|\le b_{\alpha,\beta}$
for $\forall \alpha,\,\beta$ and in this case $\phi$ is called a majorant 
series of $\psi$.
Note that if $\phi$ is a convergent power series, so is $\psi$.

Now assume \eqref{eq:converg}.
We note that there exists $\epsilon>0$ such that
\[
  |\det\bar\sigma_*(P)(0,\xi)|\ge\epsilon(\xi_1+\cdots+\xi_n)^{mr}
  \quad\text{for }\forall\xi\in[0,\infty)^n.
\]
As in the proof of \cite[Theorem~2.1]{OV1},
we can choose $C>0$, $c>0$, $M>0$ and $K\ge 1$ so that
for $\forall(\alpha,\beta)\in\mathbb N^{n+n'}$ and
$\forall \gamma\in\mathbb N^n\setminus\Sigma$
\begin{gather}
  cm|\bigl(\sigma_*(P)(0,\gamma)^{-1}\bigr)_{ij}|\le
    \prod_{j=0}^{r-1}\bigl(r|\gamma|-j\bigr)^{-1}
    ,\notag\\
  \sigma_*(P)(x,\gamma)_{ij}-\sigma_*(P)(0,\gamma)_{ij}
     \le \frac{C(x_1+\cdots+x_{n'})\prod_{j=0}^{r-1}
     \bigl(r|\gamma|-j\bigr)}{1-K(x_1+\cdots+x_{n'})}
    ,\notag\\
  p_{\alpha,\beta}(t,x)_{ij}-p_{\alpha,\beta}(0,x)_{ij}\ll
    \frac{C(t_1+\cdots+t_n)}{1-K(t_1+\cdots+t_n+x_1+\cdots+x_{n'})}
   \notag\\
  h(t,x)_i
  \ll\frac{M(t_1+\cdots+t_n)}{1-K(t_1+\cdots+t_n+x_1+\cdots+x_{n'})}.
  \notag
\end{gather}
Here $i$ and $j$ represent the indices of square matrices or vectors of 
size $m$.
Hence the power series $w(s,y)$ of $(s,y)$ satisfying
\begin{align}
 c\prod_{j=0}^{r-1}\Bigl(rs\frac{\p}{\p s}-j\Bigr)w &=
 \frac{Cmy}{1-Ky}\prod_{j=0}^{r-1}\Bigl(rs\frac{\p}{\p s}-j\Bigr) w\notag\\
 &\quad
  + \sum_{j+k\le r}
    \frac{Cm(n+n')^rs}{1-K(s+y)}\Bigl(s\frac{\p}{\p s}\Bigr)^j
    \Bigl(\frac{\p}{\p y}\Bigr)^kw\label{eq:major}\\
 &\quad
  + \frac{Ms}{1-K(s+y)},\notag\allowdisplaybreaks\\
 w(0,y) &=0\notag\\
\intertext{implies}
  \bigl(u(t,x) - \sum_{\alpha\in\mathbb N^n\setminus\overline\Sigma}
   u_\alpha(x)t^\alpha\bigr)_i
  &\ll w(t_1+\cdots+t_n,x_1+\cdots+x_{n'})
  \quad\text{for }1\le i\le m.\quad\notag
\end{align}
Put $s = z^r$.
Then \eqref{eq:major} changes into
\begin{equation}\label{eq:major2}
\begin{split}
 \Bigl(c-\frac{Cmy}{1-Ky}\Bigr)z^r\frac{\p^r w}{\p z^r} &=
    \sum_{j+k\le r}
    \frac{Cm(n+n')^rz^r}{1-K(z^r+y)}
    \Bigl(\frac zr\frac{\p}{\p z}\Bigr)^j
    \frac{\p^k w}{\p y^k}\\ 
    &\quad+\frac{Mz^r}{1-K(z^r+y)},\\
 \frac{\p^jw}{\p z^j}\Bigm|_{z=0}& = 0
 \quad\text{for}\quad j=0,\ldots,r-1.
\end{split}
\end{equation}
Since the first equation in the above is equivalent to
\[
 \Bigl(c-\frac{Cmy}{1-Ky}\Bigr)\frac{\p^r w}{\p z^r} =
     \sum_{j+k\le r}
    \frac{Cm(n+n')^r}{1-K(z^r+y)}\Bigl(\frac zr\frac{\p}{\p z}\Bigr)^j
    \frac{\p^k w}{\p y^k}+\frac{M}{1-K(z^r+y)},
\]
\eqref{eq:major2} has a unique solution of power series of $(y,z)$,  
which is assured to be analytic at the origin by Cauchy-Kowalevsky's theorem.
In fact for a sufficiently large positive number $L$, 
the solution of the ordinary differential equation
\begin{align*}
  \Bigl(c -\frac{Cmt}{1- Kt}\Bigr)\tilde w^{(r)}(t) &=
  \sum_{j+k\le r}\frac{Cm(n+n')^rL^{-k}}{1-Kt}
    \Bigl(\frac tr\frac{d}{d t}\Bigr)^j
    \tilde w^{(k)}(t)+\frac M{1-Kt},\\
  \tilde w^{(j)}(0)&=0\quad\text{for \ }j=0,\ldots,r-1
\end{align*}
with
\begin{equation*}
  \begin{cases}
    t = z+Ly,\\
    cL^r>Cm(n+n')^r
  \end{cases}
\end{equation*}
satisfies $w(z,y)\ll \tilde w(z+Ly)$.
Hence $u$ is also a convergent power series.
\end{proof}

Let $\ell$ be a non-negative integer and let $U$ be an open connected 
neighborhood of a point $z^o$ of $\mathbb C^\ell$ and let $\mathcal O_U$ 
be the space of holomorphic functions on $U$.
We denote by ${}_U\! \mathcal A_M$ and ${}_U\! \mathcal A_N$ 
the space of real analytic functions on $M$ with holomorphic parameter 
$z\in U$ and that on $N$ with holomorphic parameter $z\in U$, respectively.
Moreover we denote by ${}_U\!\hat{\mathcal A}_M$ the space of formal power 
series of $t=(t_1,\ldots,t_n)$ with coefficients in ${}_U\!A_N$.
Let ${}_U\!\mathcal D_*$ denote the ring of differential operators $P$ of the 
form
\begin{equation*}
 \begin{cases}
 P = \sum_{(\alpha,\beta)\in\mathbb N^{n+n'}}
     a_{\alpha,\beta}(t,x,z)\vartheta^\alpha\p_x^\beta,\\
 a_{\alpha,\beta} \in {}_U\!\mathcal A_N,\ 
 a_{\alpha,\beta}(0,x,z) = 0
 \quad\text{if \ }\beta>0.
 \end{cases}
\end{equation*}
Then $\sigma_*(P)(x,z,\xi) :=\sum_{\alpha}p_{\alpha,0}(0,x,z)\xi^\alpha\in{}_U\!\mathcal A_N[\xi]$.
\begin{theorem}\label{thm:ideal}
Let $P\in M(m,{}_U\!\mathcal D_*)$
and
$\lambda(z)=\bigl(\lambda_1(z),\ldots,\lambda_n(z)\bigr)
 \in \mathcal O_U^n$. 

{\rm i)}
Let $\Sigma$ be a subset of $\mathbb N^n$ such that
\begin{equation*}
  \det\bigl(\sigma_*(P)(x,z,\lambda(z)+\gamma)\bigr)\not\equiv 0
  \quad\text{for \ }\forall\gamma\in\mathbb N^n\setminus\Sigma.
\end{equation*}
Let 
$\phi(t,x,z)=\sum_{\alpha\in\mathbb N^n}
\phi_{\alpha}(x,z)t^\alpha\in{}_U\!\hat{\mathcal A}_M^m$
satisfying  $P\bigl(t^{\lambda(z)}\phi\bigr) = 0$.
Then $\phi=0$ if $\phi_\alpha=0$ for $\forall\alpha\in\Sigma$.

Hereafter in this theorem suppose $P$ satisfies
\begin{align}\label{eq:converu}
 \det \bar\sigma_*(P)(x,z,\xi)\ne0
 \quad\text{for }\forall (x,z,\xi) \in N\times U\times
 \bigl\{[0,\infty)^n\setminus\{0\}\bigr\}.
\end{align}

{\rm ii)} If $\phi(t,x,z)\in{}_U\!\hat{\mathcal A}_M^m$ satisfies
   $P\bigl(t^{\lambda(z)}\phi\bigr) = 0$,
then $\phi\in {}_U\!{\mathcal A}_M^m$.

{\rm iii)}
Fix $x^o\in N$.
Let $\Sigma$ be 
a finite subset $\Sigma$ of $\mathbb N^n$ such that
\begin{equation*}
  \det\bigl(\sigma_*(P)(x^o,z^o,\lambda(z^o)+\gamma)\bigr)\ne 0
  \quad\text{for \ }\forall\gamma\in\mathbb N^n\setminus\Sigma.
\end{equation*}
Shrinking $U$ and $N$ if necessary 
and denoting
\begin{align}
 \Sol_U(P;\lambda):&=
  \{u\,;\, ut^{-\lambda(z)}\in {}_U\!\mathcal A_M^m\text{ and }Pu
  = 0\},
\notag\\
 \begin{split}
 \Sol_U(P;\lambda)^{\Sigma}:&=
  \{\bar u=\sum_{\alpha\in\overline\Sigma}\phi_\alpha(x,z)t^{\lambda(z)+\alpha}\,;
  \,\bar ut^{-\lambda(z)}\in {}_U\!\mathcal A_M^m \text{ and }
\notag
\\ 
 &\qquad P\bar u
  \equiv 0\mod
 \sum_{\beta\in\p\Sigma}{}_U\!\mathcal A_M^mt^{\lambda(z)+\beta}\},
 \end{split}
\end{align}
we see that the natural restriction map
\begin{equation*}
 \begin{split}
 \Sol_U(P;\lambda)\, &\xrightarrow{\sim}\, \Sol_U(P;\lambda)^{\Sigma},\\
 \sum_{\alpha\in\mathbb N^n}
 \phi_\alpha(x,z)t^{\lambda(z)+\alpha}&\mapsto \sum_{\alpha\in\Sigma}
 \phi_\alpha(x,z)t^{\lambda(z)+\alpha}
\end{split}
\end{equation*}
is a bijection. 
Here in particular
\begin{equation*}
 \Sol_U(P;\lambda)^{\{0\}} 
 =\{u\in {}_U\!\mathcal A_N^m\,;\,
   \sigma_*(P)\bigl(x,z,\lambda(z)\bigr)u = 0\}.
\end{equation*}
\end{theorem}
\begin{proof}
Fix $x^o\in N$.
Expanding functions in convergent power series of $(t,x,z)$ at 
$(0,x^o,z^o)$, we will prove the lemma in a neighborhood of 
$(0,x^o,z^o)$.
Replacing $P$ and the complexification $M_{\mathbb C}$ of $M$ by
$t^{-\lambda(z)}\circ P\circ t^{\lambda(z)}$
and $M_{\mathbb C}\times U$, respectively, we can reduce this theorem 
to the previous theorem without the parameter $z$.
\end{proof}
\begin{corollary}
Retain the notation in the previous theorem.
Let $\ell=1$.
Suppose
\begin{equation*}
  \sigma_*(P)(x,z,\lambda(z))=0
  \quad\text{for \ }\forall(x,z)\in N\times U
\end{equation*}
and
\begin{equation*}
  \det\bigl(\sigma_*(P)(x^o,z,\lambda(z)+\gamma)\bigr)\ne 0
  \quad\text{for \ }\forall\gamma\in\mathbb N^n\setminus\{0\}
  \text{ and }\forall z\in U\setminus\{z^o\}.
\end{equation*}
Then there exists a non-negative integer $k$ such that the following
holds.

The previous theorem assures that for any
$\phi_0(x,z)\in{}_U\!\mathcal A_N^m$ and fixed $z\in U\setminus\{z^o\}$ 
there exists a function $u(t,x,z)$ satisfying
\[
  \begin{cases}
   P u = 0,\\
   t^{-\lambda(z)}u\in\!\mathcal A_M^m,\\
   t^{-\lambda(z)}u|_{t=0}=\phi_0(x,z).
  \end{cases}
\]
Then $t^{-\lambda(z)}z^k u(x,z)$ extends holomorphically to the point $z=z^o$.
\end{corollary}
\begin{proof}
Since the functions $\det\bigl(\sigma_*(P)(x^o,z,\lambda(z)+\gamma)\bigr)$
have finite order of zeros at $z=z^0$ for $\gamma\in\Sigma\setminus\{0\}$,
this corollary follows from the proof of Theorem~\ref{thm:sol} (cf.~\eqref{eq:zen0} for $\forall\alpha^o\in\mathbb N^n\setminus\{0\}$).
In fact it is sufficient to put $k$ the sum of these orders of zeros for 
$\gamma\in\Sigma\setminus\{0\}$. 
\end{proof}
\begin{remark}\label{rem:termD}
It follows from the proves of Theorem~\ref{thm:sol} and 
Theorem~\ref{thm:ideal} that
there exist differential operators $P_\alpha^\gamma(x,z,\p_x)$ such that
\begin{equation*}
  \phi_\alpha(x,z) = 
  \sum_{\gamma\in\Sigma}P_\alpha^\gamma(x,z,\p_x)\phi_\gamma(x,z)
  \quad\text{for }\alpha\in\mathbb N^n\setminus\Sigma
\end{equation*}
in Theorem~\ref{thm:ideal} {\rm iii)}.
\end{remark}
\begin{corollary}\label{cor:ideal}
Fix $(x^o,\lambda^o)\in N\times \mathbb C^n$
and let $V$ be a neighborhood of $\lambda^o$ in $\mathbb C^n$.
Suppose $P\in M(m,\mathcal D_*)$ satisfies \eqref{eq:converu}
and
\begin{equation*}
 \det\bigl(\sigma_*(P)(x^o,\lambda^o+\gamma\bigr)-
       \sigma_*(P)(x^o,\lambda^o\bigr)\bigr)
  \ne 0\quad\text{for \ }\forall\gamma\in\mathbb N^n\setminus\{0\}.
\end{equation*}
Then shrinking $N$, $M$  and $V$ if necessary, we have a linear bijection
\begin{align*}
  \beta_\lambda: 
   \Sol_V(P):=\{u\,;\, ut^{-\lambda}\in {}_V\! \mathcal A_M^m\text{ and }
    Pu = \sigma_*(P)(x,\lambda)u\}\ 
  \xrightarrow{\sim}& \ {}_V\! \mathcal A_N^m,\\
  u\mapsto&\ t^{-\lambda}u|_{t=0}
\end{align*}
with the coordinate $((t,x),\lambda)\in M\times V$.
In particular, we have a bijective map
\begin{align*}
  \beta_{\lambda^o}: \Sol_{\lambda^o}(P):=
  \{u\,;\, ut^{-\lambda^o}\in \mathcal A_M^m\text{ and }
    Pu = \sigma_*(P)(x,\lambda^o)u\}\ 
  \xrightarrow{\sim}& \ \mathcal A_N^m,\\
  u\mapsto&\ t^{-\lambda^o}u|_{t=0}.
\end{align*}
\end{corollary}
\begin{definition}
The map $\beta_{\lambda^o}$ of $\Sol_{\lambda^o}(P)$ is called
the \emph{boundary value map} of the solution space $\Sol_{\lambda^o}(P)$
of the differential equation $Pu=\sigma_*(P)(x,\lambda^o)u$ with respect to
the \emph{characteristic exponent} $\lambda^o$.
\end{definition}
\begin{remark}
When $n=1$, $u\in\Sol_{\lambda^o}(P)$ is called an
\emph{ideally analytic solution} of the equation
$   P u = \sigma_*(P)(x,\lambda^o)u $
in \cite{KO}.
\end{remark}
The following theorem says that $\Sol_V(P)$ and 
$\sigma_*(P)$ characterize
$P\in\mathcal D_*$.
\begin{theorem}\label{thm:sol2eq}
Let $P$ be an element of $M(m,\mathcal D_*)$
satisfying the assumptions in 
Corollary~\ref{cor:ideal}.
Let $P'\in M(m,\mathcal D_*)$ with $\sigma_*(P)=\sigma_*(P')$.
Then the condition $\Sol_V(P) = \Sol_V(P')$ implies $P=P'$.
\end{theorem}
\begin{proof}
Suppose $P\ne P'$.
Put
\[
 P - P' = \sum_{\alpha,\beta,\gamma}r_{\alpha,\beta,\gamma}
 t^\gamma\vartheta^\alpha\p_x^\beta.
\]
Then we can find $\gamma^o\in\mathbb N^{n'}\setminus\{0\}$
such that $\sum_{\alpha,\beta}r_{\alpha,\beta,\gamma^o}t^{\gamma^o}\vartheta^\alpha\p_x^\beta\ne0$ and $r_{\alpha,\beta,\gamma}=0$ if $\gamma<\gamma^o$.
For $v(x)\in\mathcal A_N^m$ the coefficients of $t^{\lambda+\gamma^o}$ in
$(P-P')\beta_{\lambda}^{-1}v(x)$ show
\begin{align*}
 0 &= \bigl(t^{-\lambda}\sum_{\alpha,\beta}r_{\alpha,\beta,\gamma^o}\vartheta^\alpha\p_x^\beta t^\lambda v(x)\bigr)|_{t=0}\\
  &= \sum_{\alpha,\beta}r_{\alpha,\beta,\gamma^o}\lambda^\alpha\p_x^\beta v(x)
 \quad\text{for \ }\forall\lambda\in V\text{ and }\forall v(x)\in\mathcal A_N^m,
\end{align*}
which means a contradiction.
\end{proof}

\section{Induced equations}\label{sec:Ind}
Retain the notation in the previous section.
Moreover we denote by ${}_U\!\widetilde{\mathcal D}_*$ the ring of holomorphic
maps of $U$ to $\widetilde{\mathcal D}_*$
for a connected open subset $U$ of $\mathbb C^\ell$.

We recall that the 
element $P$ of ${}_U\!\widetilde{\mathcal D}_*$ is characterized 
by the expression
\begin{equation}
  P = \sum_{(\alpha,\beta)\in\mathbb N^{n+n'}}
   p_{\alpha,\beta}(t,x,z)\vartheta^\alpha\p_x^\beta
\end{equation}
with $p_{\alpha,\beta}(t,x,z)\in {}_U\!\mathcal A_M$
and
\begin{equation*}
 \sigma_*(P)(x,z,\xi,\p_x) = \sum_{\alpha,\beta} p_{\alpha,\beta}(0,x,z)
 \xi^\alpha\p_x^\beta.
\end{equation*}

\begin{theorem}\label{thm:induce}
Let $P\in M(m,{}_U\!\mathcal D_*)$ satisfying the assumption 
in Theorem~\ref{thm:ideal}~iii) with $\Sigma=\{0\}$.
Suppose that $P_1,\ldots,P_p\in M(m,{}_U\!\widetilde{\mathcal D}_*)$ satisfy
\begin{equation}\label{eq:involutive}
 [P,P_i] = S_iP + \sum_{j=1}^pT_{ij}P_j
\end{equation}
with $S_i\in M(m,{}_U\!\widetilde{\mathcal D}_*)$ and 
$T_{ij}\in M(m,{}_U\!\mathcal D_*)$.
Suppose moreover $\sigma_*(T_{ij})=0$.
Then the map
\begin{equation}\label{eq:bvind}
\begin{split}
\beta_{\lambda(z)}:&\bigl\{u\,;\,t^{-\lambda(z)}u\in{}_U\!\mathcal A_M^m\text{ and }
 Pu=P_iu=0\text{ for }i=1,\ldots,p\bigr\}\\
&
\xrightarrow{\sim}\bigl\{v\in{}_U\!\mathcal A_N^m\,;\,
\begin{cases}
 \sigma_*(P)\bigl(x,z,\lambda(z)\bigr)v=0,\\
 \sigma_*(P_i)\bigl(x,z,\lambda(z),\p_x\bigr)v=0\quad(i=1,\ldots,p)
\end{cases}\bigr\},\\
&
u\mapsto 
t^{-\lambda(z)}u\bigm|_{t=0}
\end{split}
\end{equation}
is a bijection.
\end{theorem}
\begin{proof}
Since $(t^{-\lambda(z)}P_ju)|_{t=0}=\sigma_*(P_j)(x,z,\lambda(z),\p_x)t^{-\lambda(z)}u|_{t=0}$, Theorem~\ref{thm:ideal} assures that we have only to prove the 
surjectivity of the map to get the theorem.

For a given $v$ in the element of the set, we have $u\in t^{\lambda(z)} 
{}_U\!\mathcal A_M^m$ such 
that $Pu=0$ and $t^{-\lambda(z)}u|_{t=0}=v$. Then
$ PP_iu=\sum_{j=1}^p T_{ij}P_ju$, namely,
\begin{equation*}
 \begin{pmatrix}
  P - T_{11} & -T_{12} & -T_{13}   &\cdots & -T_{1p}\\
  -T_{21} &   P - T_{22} & -T_{23} & \cdots &-T_{2p}\\
  -T_{31} &   - T_{32} & P -T_{33} & \cdots &-T_{2p}\\
  \vdots  & \vdots & \vdots & \ddots & \vdots\\
  -T_{p1} & T_{p2} & T_{p3} & \cdots & P - T_{pp}
 \end{pmatrix}
 \begin{pmatrix}
  P_1u \\ P_2u \\ P_3u\\ \vdots \\ P_pu
 \end{pmatrix}
 =0.
\end{equation*}
Since $\sigma_*(T_{ij}) = 0$ and $t^{-\lambda(z)}P_ju|_{t=0} = 0$ for 
$j=1,\ldots,p$, Theorem~\ref{thm:ideal} i) assures $P_ju=0$.
\end{proof}
\begin{definition}
The system of differential equations
\begin{equation*}
 \sigma_*(P)(x,z,\lambda(z)) v = 
 \sigma_*(P_i)(x,z,\lambda(z),\p_x)v = 0\quad\text{for \ }i=1,\ldots,p
\end{equation*}
in Theorem~\ref{thm:induce} is called the system of \emph{induced equations} 
with respect to the boundary value map $\beta_{\lambda(z)}$ 
{\rm (cf.~\eqref{eq:bvind})}.
\end{definition}
\begin{remark}\label{rem:nonsplit}
{\rm i)}
Suppose $P\in M(m,{}_U\!\mathcal D_*)$ satisfies the assumption in
Theorem~\ref{thm:induce}.
Let $Q\in M(m,{}_U\!\mathcal D_*)$ such that $[P,Q]=0$ and
$\sigma_*(Q)\bigl(x,z,\lambda(z)\bigr)=0$.
Then if $u\in t^{\lambda(z)}{}_U\!\mathcal A_M^m$ satisfies $Pu=0$,
we have $Qu=0$.

{\rm ii)}
Let $p$ be the rank of an irreducible semisimple symmetric space $G/H$.
The ring of invariant differential operators on $G/H$ is isomorphic
to $\mathbb C[P_1,\ldots,P_p]$, where $P_j$ are algebraically independent
and satisfy $[P_i,P_j]= 0$ for $1\le i<j\le p$.
Under a suitable coordinate system $(t_1,\ldots,t_n,x_1,\ldots,x_{n'})$ 
of a natural realization of $G/H$ constructed by \cite{OR2}, $G/H$ is defined
by $t_1>0,\ldots,t_n>0$. 
Then $n$ is the real rank of $G/H$ and 
$P_i\in\tilde{\mathcal D}_*\setminus\mathcal D_*$ if $n < p$.
It is shown in \cite{OR2} that we can choose 
$P\in\sum_{j=1}^p\mathcal D_* P_j$ such that 
$P$, $P_1,\ldots,P_p$ satisfy the assumption in Theorem~\ref{thm:induce}.
\end{remark}
\section{Holonomic systems of differential equations 
with constant coefficients}\label{sec:HS}
In this section $\bigl(\frac{\p}{\p y_1},\ldots,\frac{\p}{\p y_n}\bigr)$
is simply denoted by $\p$.
For $\mu=(\mu_1,\ldots,\mu_n)\in\mathbb C^n$
and $y=(y_1,\ldots,y_n)\in\mathbb R^n$, we put
\begin{equation*}
  \langle\mu,y\rangle=\mu_1y_1+\cdots+\mu_ny_n.
\end{equation*}
\begin{lemma}
Let $\Hom_{\mathbb C[\p]}(\mathcal M,\mathcal N)$ denote the space of
$\mathbb C[\p]$-homomorphisms of a $\mathbb C[\p]$-module
$\mathcal M$ to a $\mathbb C[\p]$-module $\mathcal N$.
Then the space is naturally a $\mathbb C[\p]$-module.
Let $\hat{\mathcal O}$ be the space of formal power series of 
$y=(y_1,\ldots,y_n)$ and let $\mathcal O(\mathbb C^n)$ be the space of entire functions on $\mathbb C^n\ni y$.
Suppose $\mathcal M$ is a finite dimensional $\mathbb C[\p]$-module.
Then 
\begin{gather}
 \begin{split}
  \bigoplus_{\lambda\in\mathbb C^n}\Hom_{\mathbb C[\p]}\bigl(\mathcal M,
  \mathbb C[y]e^{\langle \lambda,y\rangle}\bigr)
   &\xrightarrow{\sim}
  \Hom_{\mathbb C[\p]}\bigl(\mathcal M,\bigoplus_{\lambda\in\mathbb C^n}
  \mathbb C[y]e^{\langle \lambda,y\rangle}\bigr)\\
   &\xrightarrow{\sim}
  \Hom_{\mathbb C[\p]}\bigl(\mathcal M,\mathcal O(\mathbb C^n)\bigr)\\
   &\xrightarrow{\sim}
  \Hom_{\mathbb C[\p]}\bigl(\mathcal M,\hat{\mathcal O}\bigr),
 \end{split}\label{eq:Msol}\\
  \dim \Hom_{\mathbb C[\p]}\bigl(\mathcal M,\mathcal O(\mathbb C^n)\bigr)
   = \dim\mathcal M.\label{eq:Msoldim}
\end{gather}
If $\mathcal M'$ is a quotient $\mathbb C[\p]$-module of $\mathcal M$ such that
\[
\Hom_{\mathbb C[\p]}\bigl(\mathcal M',\mathcal O(\mathbb C^n)\bigr)
  \xrightarrow{\sim}\Hom_{\mathbb C[\p]}\bigl(\mathcal M,\mathcal O(\mathbb C^n)\bigr),
\]
then $\mathcal M\xrightarrow{\sim}\mathcal M'$.
\end{lemma}
\begin{proof}
For $\mu=(\mu_1,\ldots,\mu_n)\in\mathbb C^n$, 
let $\mathfrak m_\mu$ denote the maximal ideal of $\mathbb C[\p]$ 
generated by $\frac{\p}{\p y_i}-\mu_i$ with $i=1,\ldots,n$. 
Then we have 
$\mathcal M\simeq\mathcal M_{\lambda_1}\oplus
\cdots\oplus\mathcal M_{\lambda_m}$ 
with suitable $\lambda_\nu=(\lambda_{\nu,1},\ldots,\lambda_{\nu,n})
\in\mathbb C^n$
and $\mathbb C[\p]$-modules $\mathcal M_{\lambda_\nu}$ satisfying
$\mathfrak m_{\lambda_\nu}^k\mathcal M_{\lambda_\nu}=0$ for a 
large positive integer $k$.
Hence we have only to prove the lemma for each $\mathcal M_{\lambda_\nu}$.
By the outer automorphism $\frac{\p}{\p y_i}\mapsto\frac{\p}{\p y_i}
+\lambda_{\nu,i}$ for $i=1,\ldots,n$
which corresponds to the multiplication of the functions in 
$\mathcal O(\mathbb C^n)$ or $\hat{\mathcal O}$ by 
$e^{-\langle\lambda_\nu,x\rangle}$ we may assume 
$\mathfrak m_0^k\mathcal M=0$.

Suppose $\mathfrak m_0^k\mathcal M=0$. Then
$\Hom_{\mathbb C[\p]}\bigl(\mathcal M,\mathbb C[y]\bigr)\xrightarrow{\sim}
  \Hom_{\mathbb C[\p]}\bigl(\mathcal M,\hat{\mathcal O}\bigr)$ and
\eqref{eq:Msol} is clear.
Since $\hat{\mathcal O}$ is the dual space of $\mathbb C[\p]$ by the bilinear 
form $\langle P(\p),u\rangle=P(\p)u|_{x=0}$, \eqref{eq:Msoldim} is clear. 
The last statement follows from \eqref{eq:Msoldim}.
\end{proof}
\begin{definition}\label{def:semismple}
A finite dimensional $\mathbb C[\p]$-module $\mathcal M$ is \emph{semisimple}
if
\[
 \Hom_{\mathbb C[\p]}\bigl(\mathcal M,\bigoplus_{\lambda\in\mathbb C^n}
 \mathbb Ce^{\langle\lambda,y\rangle}\bigr)
 \xrightarrow{\sim}
 \Hom_{\mathbb C[\p]}\bigl(\mathcal M,\mathcal O(\mathbb C^n)\bigr).
\]
\end{definition}
Let $U$ be a convex open subset of $\mathbb C^\ell$,
where $\ell$ is a non-negative integer, and
let ${}_U\!\mathbb C[\p]$ and ${}_U\!\mathcal O(\mathbb C^n)$ be the space
of holomorphic maps of $U$ to $\mathbb C[\p]$ and that of $U$ to 
$\mathcal O(\mathbb C^n)$, respectively.
\begin{proposition}\label{prop:cdif}
Let $r$ be a positive integer and let ${}_U\!\mathcal M$ be a 
finitely generated ${}_U\!\mathbb C[\p]$ module 
with $\dim {}_U\!\mathcal M=r$ for any fixed $z\in U$.
Assume that there exist positive integer $k$ and finite number of holomorphic
maps $\lambda_i$ of $U$ to $\mathbb C^n$ such that
$\bigl(\prod_{i\in I}\mathfrak m_{\lambda_i(z)}^k\bigr){}_U\!\mathcal M=0$
for any $z\in U$.
Here the indices $i$ run over a finite set $I$.
Then there exist ${}_U\!\mathbb C[\p]$-homomorphisms 
$u_1,\ldots,u_r$ of ${}_U\!\mathcal M$ to ${}_U\!\mathcal O(\mathbb C^n)$
such that they are linearly independent for any fixed $z\in U$.

Let $I=I_1\cup\cdots\cup I_L$ be a decomposition of $I$ such that
\begin{align*}
 \lambda_i(z)\ne\lambda_j(z)\quad{for \ }\forall z\in U
 \quad\text{if \ }i\in I_\mu\text{ and }j\in I_\nu
 \text{ \ and \ }1\le \mu<\nu\le L.
\end{align*}
Then we can choose $\{u_i\,;\,i\in I\}$ such that for each $u_i$ there 
exists $I_\nu$ satisfying
\begin{align}\label{eq:separate}
 u_i\in\Hom_{{}_U\!\mathbb C[\p]}\bigl({}_U\!\mathcal M,
\sum_{j\in I_\nu}e^{\langle\lambda_j(z),y\rangle}\mathbb C[y]\bigr)
 \quad\text{for any fixed }z\in U.
\end{align}
\end{proposition}
\begin{proof}
Let $\{v_1,\ldots,v_m\}$ be a system of generators of ${}_U\!\mathcal M$.
We identify the homomorphisms of ${}_U\!\mathcal M$ to 
${}_U\mathcal O(\mathbb C^n)$ with their image of $\{v_1,\ldots,v_m\}$
and hence $u_j(y,z)\in\mathcal {}_U\!\mathcal O(\mathbb C^n)^m$.
Note that we can find
${}_U\!\mathbb C[\p]$-homomorphisms $\tilde u_1(y,z),\ldots,\tilde u_r(y,z)$ 
of ${}_U\!\mathcal M$ to ${}_U\mathcal O(\mathbb C^n)$ if we replace 
$\mathcal O(U)$ by its quotient field.

Fix a point $z^o\in U$.
Let $\gamma(t)$ be a holomorphic map of $\{t\in\mathbb C\,;\,|t|<1\}$ to
$U$ such that $\gamma(0)=z^o$ and $\tilde u_j(y,\gamma(t))$ are holomorphic and
linearly independent for $0<|t|<1$.
Then \cite[Proposition~2.21]{OS} assures that there exist meromorphic 
functions $c_{ij}(t)$ such that the functions
$v_i(y,t)=\sum_{j=1}^rc_{ij}(t)\tilde u_j(y,\gamma(t))$ are 
holomorphic at $t=0$ and  that $v_1(y,0),\ldots,v_r(y,0)$ are linearly 
independent.
We can find $P_i\in\mathbb C[\p]^m$ such that
$\langle P_i,v_j\rangle = \delta_{ij}$ for $1\le i\le r$ and $1\le j\le r$.
Here we put $\langle(Q_1,\ldots,Q_m), (f_1,\ldots,f_m)\rangle
:=\sum_{\nu=1}^m Q_\nu(f_\nu)(0)$ for 
$Q_\nu\in\mathbb C[\p]^m$ and $f_\nu\in\mathcal O(\mathbb C^n)^m$.

Put $A(z)=\Bigl(\langle P_i, \tilde u_j\rangle\Bigr)
_{\substack{1\le i\le r\\1\le j\le r}}$, which is
a matrix of meromorphic functions on $U$ and $\det A(z)$ is 
not identically zero.
Let $\tilde c_{ij}(z)$ are meromorphic functions on $U$ such that
$\langle P_i,u_j\rangle=\delta_{ij}$ by putting
$u_i=\sum_{j=1}^r\tilde c_{ij}(z)\tilde u_j$.

Suppose $u_i(y,z)$ is not holomorphic at $z=z^0$.
Then there exist a positive integer $L$ and 
a holomorphic function $\tilde\gamma$ of 
$\{t\in\mathbb C\,;\,|t|<1\}$ to $U$ such that $\tilde\gamma(0)=z^0$ and 
the function $w(y,t):=t^Lu_i(y,\tilde \gamma(t))$ is holomorphically 
extended to the point $t=0$ and moreover $w(y,0)\ne0$. 
Then $w(y,0)$ defines a 
$\mathbb C[\p]$-homomorphism of ${}_U\!\mathcal M$ to 
$\mathcal O(\mathbb C^n)$  at $z=z^o$.
But $w(y,0), v_1(y,0),\ldots,v_r(y,0)$ are linearly independent because
$\langle P_i, w(y,0)\rangle=0$ for $i=1,\ldots,r$, which contradicts to
\eqref{eq:Msoldim}.

Hence for any $z^o\in U$
we can construct $u_1(y,z),\ldots,u_r(y,z)$ which are linearly independent 
and holomorphic in a neighborhood of $z^o\in U$.
Then the theorem follows from the theory of holomorphic functions with
several variables because $U$ is a convex open subset of $\mathbb C^\ell$.

Since we have a decomposition 
${}_U\!\mathcal M={}_U\!\mathcal M_1\oplus\cdots\oplus{}_U\!\mathcal M_L$
such that the module
$\bigl(\prod_{i\in I_\nu}\mathfrak m_{\lambda_i(z)}^k\bigr)
{}_U\!\mathcal M_\nu$ vanishes for $\nu=1,\ldots,L$, 
we can assume \eqref{eq:separate}.
\end{proof}
\begin{example}
Let $W$ be a finite reflection group on a Euclidean space $\mathbb R^n$.
Let $\mathbb C[p_1,\ldots,p_n]$ be the algebra of $W$-invariant polynomials
on $\mathbb R^n$.
For example,
$p_k(x) = \sum_{1\le i_1<\cdots<i_k\le n}x_{i_1}\cdots x_{i_k}$.
Then the system of differential equations
\[
  \mathcal M_\lambda: p_i(\p)u=p_i(\lambda)u\quad\text{for \ }i=1,\ldots,n
\]
with $\lambda\in\mathbb C^n$ is a fundamental example of a
${}_U\mathbb C[\p]$-module in Proposition~\ref{prop:cdif}.
Here $U=\mathbb C^n\ni\lambda$ and $r=\#W$.
The system is semisimple if and only if $w\lambda\ne\lambda$ for 
$\forall w\in W\setminus\{e\}$.
When $\lambda = 0$, the solutions of this system  are 
called harmonic polynomials for $W$.
In this case, an explicit construction of solutions is given by \cite{O-Asymp}
such that $u_1(\lambda,y),\ldots,u_r(\lambda,y)$ are entire functions
of $(\lambda,y)\in\mathbb C^{2n}$ and  linearly independent for
any fixed $\lambda\in\mathbb C^n$.
\end{example}
\begin{remark}
We will apply the result in this section to our original systems 
with the coordinates
$t_i=^{-y_i}$ for $i=1,\ldots,n$.
Then $\mathbb C[\p]$ and $^{\langle\lambda,y\rangle}f(y)$ change
into $\mathbb C[\vartheta]$ and $t^{-\lambda}f\bigl(-\log t_1,\ldots,-\log t_n
\bigr)$,
respectively.
\end{remark}
\section{Ideally analytic solutions for complete systems}\label{sec:AI}
In this section we will study the system of differential equations
\begin{align}\label{eq:comeq}
  \mathcal M: P_i u = 0\quad\text{for }i=0,1,\ldots,q
\end{align}
with $P_i\in M(m,{}_U\!\mathcal D_*)$.
Here $z\in U$ is a holomorphic parameter
and $U$ is a convex open subset of $\mathbb C^\ell$.
We assume that $\sigma_*(P_i)$ do not depend on $x\in N$.
We moreover assume that 
$P=P_0$ satisfies \eqref{eq:converg} and the system
\begin{align}\label{eq:indicial}
 \overline{\mathcal M}:
 \sigma_*(P_i)(z,\vartheta)\bar u=0\quad\text{for \ }i=0,1,\ldots,q,
\end{align}
which we call \emph{indicial equation},
satisfies the assumption of Proposition~\ref{prop:cdif}.
Then we call $\mathcal M$ a \emph{complete system of differential equations
with regular singularities along the set of walls} $\{N_1,\ldots,N_n\}$.

For a non-negative integer $k$ 
let $\mathbb C[\log t]_{(k)}$ denote the polynomial function of
$(\log t_1,\ldots,\log t_n)$ with degree at most $k$.
Put $\mathbb C[\log t] = \bigcup_{k=1}^\infty\mathbb C[\log t]_{(k)}$.
\begin{definition}
A solution $u(t,x,z)$ of $\mathcal M$ with the holomorphic parameter $z$ 
is called an \emph{ideally analytic solution} if 
$u(t,x,z)\in\bigoplus_{\lambda\in\mathbb C}
t^\lambda\mathbb C[\log t]\mathcal A_N^m$
for any fixed $z\in U$.
\end{definition}

First we will examine the system $\mathcal M$ 
without the holomorphic parameter $z$ or $U$ is a point.
Then let $\{\bar u_i=t^{\lambda_i}v_i(\log t)\,;\,i=1,\ldots,r\}$ be a 
basis of the solutions of \eqref{eq:indicial}.
Here $v_i(\xi)\in\mathbb C[\xi]$ and these $\lambda_i$ are called 
\emph{exponents} of the system $\mathcal M$.
We define
\begin{align*}
	\begin{cases}
	e(\bar u_i)    := \lambda_i,\\
    \deg(\bar u_i) := \deg v_i.
    \end{cases}
\end{align*}
We may assume that for any $\lambda\in\mathbb C^n$ and $k\in\mathbb N$
\begin{multline*}
  \{\bar u_i\,;\,\bigl(e(\bar u_i),\deg(\bar u_i)\bigr)=(\lambda,k)\}
  \text{ is empty }\\
  \text{or linearly independent in the space }
  t^{\lambda}\mathbb C[\log t]_{(k)}^m/t^{\lambda}\mathbb C[\log t]_{(k-1)}^m.
\end{multline*}
\begin{definition}
Let $u(t,x)$ be an ideally analytic solution of $\mathcal M$.
Then a non-zero function 
\begin{equation}\label{eq:leading}
 w(t,x)=\sum_\nu t^\lambda p_\nu(\log t)\phi_\nu(x)
\end{equation} 
with suitable $\lambda\in\mathbb C^n$, $p_\nu(\xi)\in\mathbb C[\xi]$ and 
$\phi_\nu(x)\in\mathcal A_M^m$ is called a \emph{leading term} of $u(t,x)$ if
\begin{equation*}
 u(t,x)-w(t,x)
 \in 
 \sum_{\substack{\mu\in\mathbb C^n\\ \lambda-\mu\notin\mathbb N^n}}
    t^{\mu}\mathbb C[\log t]\mathcal A_M^m
\end{equation*}
and $\lambda$ is called a \emph{leading exponent} of this leading term.
If $\{w_1(t,x),\ldots,w_k(t,x)\}$ is the complete set of the leading terms of
$u(t,x)$, we say $\sum_{i=1}^k w_i(t,x)$ the \emph{complete leading term}
of $u(t,x)$.
\end{definition}
Then we have the following theorem.
\begin{theorem}\label{thm:nopara}
The leading term \eqref{eq:leading} of an ideally analytic solution 
$u(t,x)$ of $\mathcal M$ is a solution
of \eqref{eq:indicial}.
Hence there exist $\phi_i(x)\in\mathcal A_M$ such that
\begin{align}\label{eq:idw1}
  w(t,x) = \sum_{\lambda_i=\lambda}\bar u_i(t)\phi_i(x).
\end{align}
In particular, $\lambda$ is an exponent of $\mathcal M$.

Assume
\begin{align}\label{eq:sepcond}
 \det\sigma_*(P_1)(e(\bar u_i)+\gamma)\ne 0
 \quad\text{for \ }\gamma\in\mathbb N^n\setminus\{0\}.
\end{align}
Then for any $\phi(x)\in\mathcal A_N$ there is a unique solution
of $\mathcal M$ in the space $t^{e(\bar u_i)}\mathbb C[\log t]\mathcal A_M^m$
whose leading term equals $\phi(x)\bar u_i$.
Denoting the solution by $T_{\bar u_i}(\phi)$, we have the following 
bijective isomorphism if \eqref{eq:sepcond} is valid for $1\le i\le r$.
\begin{align}\label{eq:allsol}
 \mathcal A_N^r\xrightarrow{\sim}\{\text{ideally analytic solutions of }\mathcal M\},\ 
 (\phi_i)\mapsto\sum_{i=1}^rT_{\bar u_i}(\phi_i).
\end{align}
\end{theorem}
\begin{proof}
Examining the equation $Pu(t,x)=0$ modulo
$\sum_{\substack{\mu\in\mathbb C^n\\ \lambda-\mu\notin\mathbb N^n}}
t^\mu\mathbb C[\log t]\mathcal A_N^m$,
we have $\sigma_*(P)(\vartheta)w(t,x)=0$ and thus \eqref{eq:idw1}.

Put $\lambda=e(\bar u_i)$.
First suppose $\deg(\bar u_i)=0$.
Then under the condition \eqref{eq:sepcond},
Theorem~\ref{thm:sol} assures the unique existence of 
$\tilde\phi(t,x)\in\mathcal A_N^m$ such that
$P_1t^\lambda\tilde\phi(t,x)=0$ and 
$t^{e(\bar u_i)}\tilde\phi(0,x)=\phi(x)u_i(t)$
and moreover Theorem~\ref{thm:induce} assures 
$P_jt^\lambda\tilde\phi(t,x)=0$.
If there exists another solution 
$\tilde u\in t^\lambda\mathbb C[\log t]\mathcal A_N^m$ of 
$\mathcal M$ with the same property, 
the leading exponent $\lambda'$  of $u - \tilde u$
satisfies $\lambda'-e(\bar u_i)\in\mathbb N^n\setminus\{0\}$,
which contradicts to \eqref{eq:sepcond}.
Thus we have proved the required uniqueness of the solution.

Next suppose $u_i=t^\lambda v_i(\log t)$ with $\deg v_i>0$.
Let $V$ be a vector space spanned by the components of elements
of $\mathbb C[\p_\xi]v_i(\xi)$ and let $\{w_1(\xi),\ldots,w_q(\xi)\}$ be 
a basis of $V$.
Here we may assume $\mathbb C[\p_\xi]w_k\in\sum_{\nu=1}^k\mathbb C w_\nu$
for $k=1,\ldots,q$.
Let $\hat u$ be the vector of size $qm$ with components
$\hat u_\nu w_\nu(\log t)$ with $\hat u_\nu\in t^\lambda \mathcal A_N^m$
for $\nu=1,\ldots,q$.
Then the system $\mathcal M$ is replaced by a system $\hat{\mathcal M}$
with an unknown function $\hat u$ where $P_i$ are replaced by suitable
$\hat P_i\in M(qm,\mathcal D_*)$, respectively.
We note that $\hat{\mathcal M}$ also satisfies the assumption of the theorem
because $\det\bigl(\sigma_*(\hat P_i)\bigr)
=\det\bigl(\sigma_*(P_i)\bigr)^q$.
Thus we may only consider the solutions with components in 
$t^\lambda\mathcal A_M$.

For example, if $n=n'=1$ and $P=(\vartheta-\lambda)^2+t^2\p_x^2$,
the solution of the equation $Pu=0$ in the space
$t^\lambda\mathcal A_M\oplus(t^\lambda \log t)\mathcal A_M$ 
corresponds to the solution of
\[
   \begin{pmatrix}
     (\vartheta-\lambda)^2+t^2\p_x^2 & 2(\vartheta -\lambda)\\
       & (\vartheta-\lambda)^2+t^2\p_x^2
   \end{pmatrix}
   \begin{pmatrix}
    u_1\\ u_2
   \end{pmatrix}
   = 0
\]
in the space $t^\lambda\mathcal A_M^2$ by the correspondence
$u=u_1+u_2\log t$.

To complete the proof of the system we have only to prove that the map
\eqref{eq:allsol} is surjective.
Let $u$ be any ideally analytic solution of $\mathcal M$.
Then any leading exponent of $u$ is an exponent of the system $\mathcal M$ and
therefore we define $\phi_i(x)$ by \eqref{eq:idw1} if $e(\bar u_i)$ is a 
leading exponent of $u$ and by $0$ otherwise.
Then if $u\ne\sum_iT_{\bar u_i}(\psi)$,
any leading exponent of $u-\sum_i T_{\bar u_i}(\phi)$ is not in the set 
$\{e(\bar u_i)\}$, which contradicts the first claim in the theorem.
\end{proof}

We will return to the case when \eqref{eq:comeq} is the complete system which 
has a holomorphic parameter $z\in U\subset\mathbb C^\ell$.

First assume that $\overline{\mathcal M}$ is semisimple for any $z\in U$
(cf.~Definition~\ref{def:semismple}) and that 
the indicial equation $\overline{\mathcal  M}$ satisfies 
the assumption in Proposition~\ref{prop:cdif} by putting $t_i=e^{-x_i}$ 
for $1\le i\le n$.
Then the proof of the previous theorem implies the following.
\begin{proposition}
Assume that $\overline{\mathcal M}$ is semisimple for any $z\in U$.
Let $\{\bar u_i(x,z)=t^{\lambda_i(z)}f_i(z)\,;\,i=1,\ldots,r\}$
be a basis of the solutions of \eqref{eq:indicial}
for any $z\in U$.  Here $f_i(z)\in\mathcal O(U)^m$.
Assume \eqref{eq:sepcond} for any $z\in U$.
Then $T_{\bar u_i}(\phi)$ is holomorphic for $z\in U$
under the notation in Theorem~\ref{thm:nopara}.
\end{proposition}

To examine the case without the assumption in this proposition,
we study a generic holomorphic curve $t\mapsto z(t)$ through 
the point $z^o\in U$ where the assumption breaks.
Hence we restrict the case when $\ell=1$.

Suppose $\ell=1$ and fix $z^o\in U$.
For simplicity we put $z^o=0$.
Assume that $\mathcal M$ is semisimple (cf.~Definition~\ref{def:semismple})
for any fixed $z\in U\setminus\{0\}$.
We will shrink $U$ if necessary hereafter until the end of the
following theorem.
Let $\{\bar u_1,\ldots,\bar u_r\}$ be a basis of the solutions
of the indicial equation for $\forall z\in U\setminus\{0\}$,
where $\bar u_i$ are 
\begin{equation*}
  \bar u_i(t,z) = t^{\lambda_i(z)}f_i(z)\quad\text{for \ }i=1,\ldots,r
\end{equation*}
with suitable $f_i\in\mathcal O(U)^m$.
Then Proposition~\ref{prop:cdif} assures that
there exist meromorphic functions $c_{ij}(z)$ such that by
denoting
\begin{equation*}
  \bar w_i(t,z) = \sum_{j=1}^r c_{ij}(z) \bar u_j(t,z),
\end{equation*}
$\{\bar w_1,\ldots,\bar w_r\}$ is a basis of the solutions
of the indicial equation for $\forall z\in U$ and $\bar w_j(t,z)$ are 
holomorphic function of $(\log t,z)\in \mathbb C^n\times U$.
By virtue of \eqref{eq:separate},
we may assume $c_{ij}(z)=0$ if $\lambda_i(0)\ne\lambda_j(0)$.

Then we have the following theorem which is the main purpose of this note.
\begin{theorem}\label{thm:main}
Under the notation above, 
there exist differential operators $R_{ij}(x,z,\p_x)$ such that
for any $\phi(x,z)\in{}_U\!\mathcal A_M^m$,
$\sum_{i=1}^r T_{\bar u_i}\bigl(R_{ij}(x,z,\p_x)\phi(x,z)\bigr)$ 
is a holomorphic function of $z\in U$ and 
an ideally analytic solution of $\mathcal M$ with the complete leading term
$\phi(x)\bar w_i(t,z)$ for any fixed $z\in U$.
Moreover the map
\begin{equation*}
\begin{split}
 \mathcal A_N^r&\xrightarrow{\sim}
 \{\text{ideally analytic solutions of }\mathcal M\},\\
 \bigl(\phi_i(x)\bigr)&\mapsto 
 \sum_{i,j}T_{\bar u_i}\bigl(R_{ij}(x,z,\p_x)\phi_j(x)\bigr)
\end{split}
\end{equation*}
holomorphically depends on $z\in U$ and it is bijective for any $z\in U$.
Here $R_{ij}(x,z,\p_x)$ are holomorphic functions of $z\in U\setminus\{0\}$ 
valued in the space of differential operators on $N$ and may have at most 
poles at $z=0$ and moreover
\begin{equation*}
  R_{ij}(x,z,\p_x) =
  \begin{cases}
   0			  &\text{if \ }\lambda_i(0)-\lambda_j(0)\notin\mathbb N^n,\\
   c_{ij}(z)      &\text{if \ }\lambda_i(0)=\lambda_j(0).
  \end{cases}
\end{equation*}
\end{theorem}
\begin{proof}
We will inductively construct $R_{ij}(x,z,\p_x)$ according to the number
$L(\lambda_j)=\sum_{\nu=1}^n\Re\lambda_{j,\nu}(0)$.
Here $\lambda_j=(\lambda_{j,1},\ldots,\lambda_{j,n})$
and $\Re\zeta$ denotes the real part of $\zeta\in\mathbb C$.

Fix a positive integer $k$ with $k\le r$.
By the hypothesis of the induction we may assume
that $R_{ij}$ have been constructed if $L(\lambda_j)>L(\lambda_k)$.
Put $R_{jk}^{(0)}=c_{ik}(z)$.  
We inductively define $R_{ik}^{(\nu)}$ for $\nu=0,1,\ldots$ as follows. 
Put
\begin{equation*}
 \sum_{i=1}^rT_{\bar u_i}\bigl(R_{ik}^{(\nu)}\phi(x,z)\bigr)
 = z^{-n_\nu}\phi^{(\nu)}_{n_\nu}(t,x)+\cdots + z^{-1}\phi^{(\nu)}_1(t,x)
   + \phi^{(\nu)}_0(t,x,z)
\end{equation*}
with $\phi^{(\nu)}_0(t,x,z)\in {}_U\!\mathcal A_M$.
Suppose $n_\nu>0$.
By the analytic continuation of $z^{n_\nu}\sum_i
T_{\bar u_i}\bigl(R_{ik}^{(\nu)}\phi(x,z)\bigr)$, 
it is clear that 
$\phi^{(\nu)}_{n_\nu}(t,x)$ is a solution of $\mathcal M$ at $z=0$.
Any leading exponent $\mu$ of $\phi^{(\nu)}_{n_\nu}(t,x)$ satisfies
$\mu-\lambda_k(0)\in\mathbb N^n\setminus\{0\}$ and hence the complete leading
term of $\phi^{(\nu)}_{n_\nu}(t,x)$ is
\[
 \sum_{\lambda_j(0)\in\lambda_k(0)
  +\bigl(\mathbb N^n\setminus\{0\}\bigr)}
  \psi^{(\nu)}_j(x)\bar w_j(t,0).
\]
Note that $\psi^{(\nu)}_j(x)=P^{(\nu)}_j(x,\p_x)\phi(x)$
for some differential operators which do not depend on $\phi(x)$.
Put $P^{(\nu)}_j(x,\p_x)=0$ if $\lambda_j(0)-\lambda_k(0)\notin
  \mathbb N^n\setminus\{0\}$.
Hence 
\begin{equation}\label{eq:polesInd} 
\sum_{i=1}^rT_{\bar u_i}\bigl(R_{ik}^{(\nu)}(x,z,\p_x)\phi(x)\bigr)
  -\sum_{i=1}^r\sum_{j=1}^r
   z^{-n_\nu}T_{\bar u_i}
\bigl(R_{ij}(x,z,\p_x)P^{(\nu)}_j(x,\p_x)\phi(x)\bigr)
\end{equation}
has a pole of order less than $n_k$.
Defining
\[
  R^{(\nu+1)}_{ik}(x,z,\p_x) = R_{ik}^{(\nu)}(x,z,\p_x)
   -\sum_{j=0}^r
   z^{-n_\nu}R_{ij}(x,z,\p_x)P_j^{(\nu)}(x,\p_x)
\]
inductively, we have $R_{ij}(x,z,\p_x)=R_{ij}^{(\nu)}(x,z,\p_x)$ for 
certain $\nu$ such that the left hand size of \eqref{eq:polesInd} is
holomorphic at $z=0$.
\end{proof}
\begin{remark}
Let $P_i\in\mathcal D_*$ for $i=1,\ldots,n$ satisfies
\begin{equation*}
 \begin{cases}
 [P_i,P_j] = \sum_{\nu=1}^n  R_{ij\nu} P_\nu
  \quad\text{for \ }1\le i\le j\le n,\\
 \sigma_*(P_i)\text{ do not depend on }x\in N,\\
 \{\xi\in\mathbb C^n\,;\,\bar\sigma_*(P_1)(\xi)
  =\cdots=\bar\sigma_*(P_n)(\xi)=0\}=\{0\}
 \end{cases}
\end{equation*}
with some $R_{ij\nu}\in\mathcal D_*$ satisfying $\sigma_*(R_{ij\nu})=0$.
Then for a suitable  positive integer $L$ there exist 
$R_i\in\mathbb C[\vartheta]$ such that
\[
 \begin{cases}
 \ord P_i +\ord R_i = 2L,\\
 \sigma_*(P_0) = \xi_1^{2L}+\cdots+\xi_n^{2L}
 \end{cases}
\]
by putting
\[
 P_0 = \sum_{i=1}^n R_iP_i.
\]
Then $\{P_0,\ldots,P_n\}$ satisfies \eqref{eq:involutive} with
$S=0$ and $\sigma_*(T_{ij})=0$ because
\[
 [P_0,P_j] = \sum_{i=1}^n \bigl([P_0,R_j]P_i
 +\sum_{\nu=1}^nR_iR_{ij\nu}P_\nu\bigr)
\]
and $\sigma_*([P_0,R_i])=\sigma_*(R_iR_{ij\nu})=0$.

In this case let $\lambda^o$ be an exponent of the system
$P_iu=0$ $(1\le i\le n)$. 
Then for a suitable $\rho\in\mathbb C^n$ and a positive integer $k$, 
the system 
\[
  {}_U\!\mathcal M:\,\bigl(P_i-\sigma_*(P)(\lambda^o+\rho z^k)\bigr)u
 =\sum R_i\bigl(P_i-\sigma_*(P)(\lambda^o+\rho z^k)\bigr)u=0
\]
satisfies the assumption of Theorem~\ref{thm:main} 
for $U=\{z\in\mathbb C\,;\,|z|<1\}$ by changing the lower order terms 
of $R_i$ if necessary.
Hence we can analyze the ideally analytic solutions of $\mathcal M$
by the analytic continuation of the parameter $z$ to the origin.
\end{remark}
\begin{theorem}\label{thm:Hideal}
Retain the notation and the assumption in Theorem~\ref{thm:main}.
Let $r'$ be the dimension of the finitely generated 
$\mathbb C[\vartheta]$-module 
\[
\bar{\mathcal M}^o := \sum_{j=1}^m\mathbb C[\vartheta]u_j\Bigm/
 \sum_{i=0}^q\sum_{k=1}^m\mathbb C[\vartheta]\sum_{j=1}^m
 \bar\sigma_*(P_i)_{kj}(z^o,\vartheta)u_j.
\]
Suppose $n'=0$ and $r'\le r$.
Then $r'=r$ and any solution of $\mathcal M$ defined on a small connected
neighborhood of $(t^o,x^o)\in M$ with $z=z^o$ and
$0<|t^o_j|\ll 1$ for $j=1,\ldots,n$ is an ideally analytic solution given
in Theorem~\ref{thm:main}.
In particular the dimension of space of the solutions equals $r$.
\begin{proof}
Let $w_\nu$ for $\nu=1,\ldots,r'$ be elements of 
$\sum_{j=1}^m\mathbb C[\vartheta]u_j$ whose residue classes form 
a basis of $\bar{\mathcal M}^o$.
Fix $z=z^o$.
Then in a neighborhood of $(0,x^o)$
\[
  \sum_{j=1}^r \mathcal A_M[\vartheta]u_j=
  \sum_{\nu=1}^{r'}\mathcal A_M w_\nu
 + \sum_{i=0}^q\sum_{k=1}^m\mathcal A_M[\vartheta]\sum_{j=1}^m
  (P_i)_{kj}u_j.
\]
Let $w$ be a column vector of size $r'$ with components $w_\nu$.
Then the system $\mathcal M$ implies
\[
    \mathcal N\,:\ \vartheta_j w = Q_j(t) w\quad\text{for \ }j=1,\ldots,n
\]
with suitable $Q_j\in M(r',\mathcal A_M)$.
Then any solution $w(t)$ of $\mathcal N$ on a neighborhood of $(t^o,x^o)$ 
is analytic and $w=0$ if $w(t^o)=0$.
Hence the dimension of the space of solutions of $\mathcal N$ is smaller
than or equals to $r'$.
But we have constructed $r$ linearly independent solutions in 
Theorem~\ref{thm:main}.
Hence we have this theorem.
\end{proof}
\end{theorem}
\begin{remark}
Retain the notation in Theorem~\ref{thm:Hideal}.
Suppose $q=n-1$, $[P_i,P_j]=0$ for $0\le i<j\le q$, 
$\bar\sigma_*(P_i)$ are diagonal matrices and 
\[
 \{\xi\in\mathbb C^n\,;\,\bar\sigma_*(P_i)(\xi)=0\quad\text{for \ }i=0,\ldots,q
 \} = \{0\}.
\]
Then $r'=r$ and $r'=m\prod_{i=0}^q\ord P_i$.
\end{remark}
\section{Examples related to $SL(n,\mathbb R)$}\label{sec:SL}
For a connected real reductive Lie group $G$ and an open subgroup $H$ 
of the fixed point group of an involutive automorphism 
$\sigma$ of $G$,  the homogeneous space $G/H$ is called a 
\emph{reductive symmetric homogeneous space}.
Then in a suitable realization $\widetilde X$ of $G/H$ constructed by 
\cite{OR2}, the system of differential equations that defines the 
simultaneous eigenspace of the elements of the ring $\mathbb D(G/H)$ of 
the invariant differential operators on $G/H$ has regular singularities 
along the boundaries of $G/H$ in this realization.
It is an important problem to study the eigenspace.
For example, see \cite{K--} in the cases of Riemannian symmetric spaces.

Note that the Lie group $G$ is identified with a symmetric homogeneous 
space of $G\times G$ with respect to the involutive automorphism $\sigma$ 
of $G$ defined by $\sigma(g_1,g_2)=(g_2,g_1)$ for 
$(g_1,g_2)\in G_1\times G_2$ and that any irreducible admissible 
representation of $G$ can be realized in an eigenspace of $\mathbb D(G)$.

In this section we will consider differential equations related to the 
Lie group $G=SL(n,\mathbb R)$, which give examples of the differential 
equations we study in this note.
The element of the Lie algebra $\mathfrak{sl}(n,\mathbb R)$ of $G$ is 
identified with that of $M(n,\mathbb R)$ whose trace equals $0$.
Let $E_{ij}$ be the fundamental matrix unit whose $(i,j)$-component equals
1 and the other components are 0.
Then $\mathfrak{sl}(n,\mathbb R)$ is spanned by the elements
$\tilde E_{ij}=E_{ij}-\frac{\delta_{ij}}n(E_{11}+\cdots+E_{nn})$ with 
$1\le i\le j\le n$.
For simplicity we put $\tilde E_i=\tilde E_{ii}$.

We identify $\mathfrak{sl}(n,\mathbb R)$ with the space of right invariant
vector field on $G$ by
\begin{equation*}
 (Xf)(g) = \frac{d}{dt}f(ge^{tX})\Bigm|_{t=0}\quad\text{for }
X\in\mathfrak{sl}(n,\mathbb R),\ f\in C^\infty(G)\text{ and }g\in G.
\end{equation*}
Here we note that
\begin{equation*}
 (E_{pq}f)\bigl((x_{ij})\bigr)
 := \frac{d}{dt}f\bigl((x_{ij})e^{tE_{pq}}\bigr)\Bigm|_{t=0}
  = \Bigl(\sum_{\nu=1}^n x_{\nu p}\frac{\p f}{\p x_{\nu q}}\Bigr)
    \bigl((x_{ij})\bigr)
\end{equation*}
for $g\in C^\infty\bigl(GL(n,\mathbb R)\bigr)$ and
$\bigl(x_{ij}\bigr)_{\substack{1\le i\le n\\1\le j\le n}}\in GL(n,\mathbb R)$
because $(i,j)$-component of 
$\bigl(x_{ij}\bigr)_{\substack{1\le i\le n\\1\le j\le n}}E_{pq}$ equals
$x_{ip}\delta_{q j}$.

We first review, by examples, that the invariant differential operators of
the Riemannian symmetric space $G/K$ has regular singularities along the 
boundaries of the space in the realization constructed in \cite{OR1}.
By the Iwasawa decomposition $G=\bar NAK$ with
\begin{align}
 K &= SO(n) = \{g\in SL(n,\mathbb R)\,;\,{}^tgg=I_n\},\notag
\allowdisplaybreaks\\
 A &=\left\{a=\begin{pmatrix}
      a_1 &  & \\
          &\ddots & \\
          &  & a_n
     \end{pmatrix}\,;\,a_j>0\text{ \ for }1\le j\le n
     \text{ and }a_1\cdots a_n=1
     \right\},\label{eq:a}
\allowdisplaybreaks\\
 \bar N &=\left\{
          \begin{pmatrix}
           1 \\
           x_{21} & 1\\
           \vdots & \vdots & \ddots &\\
           x_{n1} & x_{n2}& \cdots & 1
          \end{pmatrix}\,;\,x_{ij}\in\mathbb R\text{ \ for }1\le j<i\le n
          \right\},\notag\\
  t_j :&= \frac{a_{j+1}}{a_j}\qquad\text{for }j=1,\ldots,n-1,\notag
\end{align}
the Riemannian symmetric space $G/K$ is identified with the product manifold
$\bar N\times A$ with the coordinate 
$(t_k,x_{ij})\in (0,\infty)^{n-1}\times\mathbb R^{\frac{n(n-1)}2}$.
Then the Lie algebra of the solvable group of $\bar NA$ is spanned by the
elements
\begin{align*}
 E_{ij} &=  \Bigl(\prod_{\nu=j}^{i-1}t_\nu\Bigr)
            \Bigl(\frac{\p}{\p x_{i j}}
             +\sum_{\nu=i+1}^nx_{\nu i}\frac{\p}{\p x_{\nu j}}\Bigr)
             \quad\text{for }1\le j<i\le n,\\
 \tilde E_{ij} &= E_{ij} - \frac{\delta_{ij}}n(E_{11}+\cdots+E_{nn})
            \quad\text{for }1\le i\le n\text{ and } 1\le j\le n,\\
 E_i :&= \tilde E_{ii} = 
 \vartheta_{i-1}-\vartheta_i\quad\text{for }1\le i\le n,\quad
           \vartheta_0 = \vartheta_{n+1}=0.
\end{align*}
The coordinate $(t_k,x_{ij})\in\mathbb R^{\frac{(n+2)(n-1)}2}$ can 
be used for local coordinate of the realization of $G/K$.

Let $U(\mathfrak g)$ be the universal enveloping algebra of the 
complexification $\mathfrak g$ of the Lie algebra of $G$.
Then if $G=SL(n,\mathbb R)$, the ring $\mathbb D(G/K)$ is naturally
isomorphic to the center $U(\mathfrak g)^G$ of $U(\mathfrak g)$ and
$U(\mathfrak g)^G$ is generated by the elements $L_2,\ldots,L_n$ which 
are given by
\begin{equation*}
\det \bigl(\tilde E_{ij}+(\tfrac{n+1}2-i-\lambda)\delta_{ij}\bigr)
=L_n - L_{n-1}\lambda + \cdots + (-1)^n\lambda^n
\end{equation*}
for $\lambda\in\mathbb C$ (cf.\ \cite{Ca}).
Here
$\det(A_{ij}) = \sum_{\sigma\in\mathfrak S_n}\sign(\sigma)A_{\sigma(1)1}\cdots A_{\sigma(n)n}$
and $U(\mathfrak g)^G$ is generated by the algebraically independent 
$(n-1)$-elements which are the coefficients of $\lambda^k$
for $k=0,1,\ldots,n-2$.

Let $\mathfrak k$ be a Lie algebra of $SO(n)$, which is generated
by the elements $E_{ij}-E_{ji}$ for $1\le i<j\le n$.

Since
\begin{align*}
 \Delta_2&=\det\begin{pmatrix} E_1 +\frac12& E_{12}\\ E_{21} & E_2-\frac12
 \end{pmatrix}
 = (E_1+\tfrac12)(E_2-\tfrac12) - E_{21}E_{12}\\
 &\equiv (E_1+\tfrac12)(E_2-\tfrac12) - E_{21}^2
   \mod U(\mathfrak g)\mathfrak k\notag\\
 &= -(\vartheta-\tfrac12)^2-t^2\p_x^2
 = -t^2(\p_t^2+\p_x^2)-\tfrac14\quad\text{with \ }\vartheta=t\tfrac \p{\p t},
 \notag
\end{align*}
we see that
$\mathbb D\bigl(SL(2,\mathbb R)/SO(2)\bigr)
=\mathbb C[t^2\bigl(\frac{\p^2}{\p t^2}+\frac{\p^2}{\p x^2}\bigr)]$.
Here $SL(2,\mathbb R)/SO(2)$ is realized in the upper half plane 
$\{x+i t\,;\,(t,x)\in(0,\infty)\times\mathbb R\}$ and
$\Delta_2$ has regular singularities along the real axis.
On the other hand, the explicit form of the vector field $L_X$ 
defined by the translation $e^{-sX}\cdot p$ for $s\in\mathbb R$, 
$X\in\mathfrak g$ and $p\in SL(2,\mathbb R)/SO(2)$ is given by
\[
 L_{E_{21}}=-\p_x,\ L_{E_1}=\vartheta+x\p_x,\ 
 L_{E_{12}}=2x\vartheta - (t^2-x^2)\p_x.
\]

When $G=SL(3,\mathbb R)$, we have
\begin{align*}
&\det
\begin{pmatrix}
 E_1+ 1 - \lambda & E_{12} &  E_{13}\\
 E_{21} & E_2-\lambda & E_{23}\\
 E_{31} & E_{32} & E_3 -1 -\lambda 
\end{pmatrix}
 = (E_1+ 1 - \lambda)(E_2-\lambda)(E_3 -1 -\lambda)\\
 &\quad
 + E_{21}E_{32}E_{13}+E_{31}E_{12}E_{23}
 - (E_{11}+1-\lambda)E_{32}E_{23} - E_{21}E_{12}(E_3-1-\lambda)\\
 &\quad-E_{31}(E_2-\lambda)E_{13}
 = \Delta_3-\Delta_2\lambda-\lambda^3
\end{align*}
with
\begin{align*}
 \Delta_3 &= (E_1+ 1)E_2(E_3 -1)
 + E_{21}E_{32}E_{13}+E_{31}E_{12}E_{23}\\
 &\quad
 - (E_1+1)E_{32}E_{23} - E_{21}E_{12}(E_3-1)
 - E_{31}E_2E_{13}\\
 &\equiv
  (E_1+ 1)E_2(E_3 -1)
 - (E_1+1)E_{32}^2 - (E_3-1)E_{21}^2 - (E_2-1)E_{31}^2\\
 &\qquad
  +2E_{21}E_{32}E_{31}
 \mod U(\mathfrak g)\mathfrak k\allowdisplaybreaks\\
 &=
  -(\vartheta_1-1)(\vartheta_1-\vartheta_2)(\vartheta_2-1)
   +2t_1^2t_2^2(\p_x+y\p_z)\p_y\p_z\\
 &\qquad
  +(\vartheta_1-1)t_2^2\p_y^2
  -(\vartheta_1-\vartheta_2-1)t_1^2t_2^2\p_z^2
  -(\vartheta_2-1)t_1^2(\p_x+y\p_z)^2,
\allowdisplaybreaks\\ 
 \Delta_2 &= E_2(E_3-1) + (E_1+1)(E_3-1) + (E_1+1)E_2
 \\&\quad 
 -E_{32}E_{23}-E_{21}E_{12}-E_{31}E_{13}
\allowdisplaybreaks\\
 &\equiv
   E_2(E_3-1) + (E_1+1)(E_3-1) + (E_1+1)E_2\\
 &\quad - E_{32}^2-E_{21}^2-E_{31}^2
   \mod U(\mathfrak g)\mathfrak k
\allowdisplaybreaks\\
 &= -(\vartheta_1-1)^2+(\vartheta_1-1)(\vartheta_2-1) -(\vartheta_2-1)^2\\
 &\quad
   -t_2^2\p_y^2-t_1^2t_2^2\p_z^2-t_1^2(\p_x+y\p_z)^2,\\
 x&=x_{21},\ y=x_{32} \text{ and }z = x_{31}.
\end{align*}
Then $\mathbb D\bigl(SL(3,\mathbb R)/SO(3)\bigr)
 =\mathbb C[\bar\Delta_3,\bar\Delta_2]$,
where
$\bar\Delta_3$ and $\bar\Delta_2$ are the last expressions of $\Delta_3$ 
and 
$\Delta_2$ in the above, respectively.
This expression of invariant differential operators on 
$SL(3,\mathbb R)/SO(3)$ is given by \cite{OV0} to obtain the Poisson 
integral representation of any simultaneous eigenfunction of the operators
on the space, where such representation is first obtained in the space with 
the rank larger than one. 
In fact $4\Delta_2$ and $8\Delta_2+8\Delta_3$ are explicitly written there 
under the coordinate $(s,t,u,v,w)$ with $(s,t,u,v,w)=(t_2^2,t_1^2,x,y,z)$, 
which corresponds to a local coordinate system in the realization given 
in \cite{OS}.

When $G=SL(n,\mathbb R)$ the second order element $L_2$ of 
$U(\mathfrak g)^G$ is
\begin{align*}
 L_2 &=
  \sum_{1\le i<j\le n}\bigl((E_i+\frac{n+1}2-i)(E_j+\frac{n+1}2-j)
   -E_{ji}E_{ij}\bigr)\\
  &\equiv  \sum_{1\le i<j\le n}(\tilde\vartheta_{i-1}-\tilde\vartheta_i)
       (\tilde\vartheta_{j-1}-\tilde\vartheta_j)\\
  &\quad-
    \sum_{1\le i<j\le n}\Bigl(\prod_{\nu=i}^{j-1}t_\nu^2\Bigr)
            \Bigl(\frac{\p}{\p x_{j i}}
             +\sum_{\nu=j+1}^nx_{\nu j}\frac{\p}{\p x_{\nu i}}\Bigr)^2
   \mod U(\mathfrak g)\mathfrak k,
 \notag\\
  \tilde\vartheta_i &= \vartheta_i-\tfrac{i(n-i)}2
\end{align*}
and $\mathbb D(G/K)=\mathbb C[\bar L_2,\ldots,\bar L_n]$ satisfying
\begin{equation*}
\begin{split}
 \sigma_*(\bar L_k)
  &= \sum_{1\le i_1<i_2<\cdots<i_k\le n}
    (\tilde\xi_{i_1-1}-\tilde\xi_{i_1})
    (\tilde\xi_{i_2-1}-\tilde\xi_{i_2})
    \cdots
    (\tilde\xi_{i_k-1}-\tilde\xi_{i_k}),\\
 \tilde\xi_i &= \xi_i-\tfrac{i(n-i)}2
\end{split}
\end{equation*}
for $k=2,\ldots,n$.

We will examine more examples.
For $a$ in \eqref{eq:a} we have
\begin{align*}
  \Ad(a^{-1})E_{ij} &:= aE_{ij}a^{-1} = a_i^{-1}a_jE_{ij} = t_{ij}E_{ij},\\
  t_{ij} &= a_i^{-1}a_j = 
    \begin{cases}
    t_it_{i+1}\cdots t_{j-1} &\text{if }i\le j,\\
    t_j^{-1}t_{j+1}^{-1}\cdots t_{i-1}^{-1} &\text{if }i>j,
  \end{cases}\\
  U(\mathfrak g)\mathfrak k &=\sum_{1\le i<j\le n}U(\mathfrak g)(E_{ij}-E_{ji}).\end{align*}
Hence
\begin{align*}
 &\Ad(a^{-1})(E_{ij}-E_{ji})^2+(E_{ij}-E_{ji})^2\\ 
 &\qquad\qquad
 -(t_{ij}+t_{ij}^{-1})\Ad(a^{-1})(E_{ij}-E_{ji})\cdot(E_{ij}-E_{ji})\\
 &\qquad=(t_{ij}^2-1)E_{ij}E_{ji}+(t_{ij}^{-2}-1)E_{ji}E_{ij}\\
 &\qquad=(t_{ij}-t_{ij}^{-1})^2E_{ji}E_{ij}+(t_{ij}^2-1)(E_{ii}-E_{jj}),\\
 &\Ad(a^{-1})(E_{ij}-E_{ji})\cdot E_{ij}-t_{ij}E_{ij}^2=
 -t_{ij}^{-1}E_{ji}E_{ij}.
\end{align*}
Thus we have
\begin{align}
 E_{ji}E_{ij} &=\frac{t_{ij}^2}{1-t_{ij}^2}(E_{ii}-E_{jj})
   \label{eq:Ktype}\\
   &\quad
   +\frac{t_{ij}^2}{(1-t_{ij}^2)^2}
            \bigl(\Ad(a^{-1})(E_{ij}-E_{ji})^2+(E_{ij}-E_{ji})^2\bigr)\notag\\
   &\quad
 -\frac{t_{ij}(1+t_{ij}^2)}{(1-t_{ij}^2)^2}\Ad(a^{-1})(E_{ij}-E_{ji})\cdot(E_{ij}-E_{ji})
  \notag\\
 &=t_{ij}^2E_{ij}^2 - t_{ij}\Ad(a^{-1})(E_{ij}-E_{ji})\cdot E_{ij}.
 \label{eq:Ntype}
\end{align}

Let $(\varpi, V_\varpi)$ be a finite dimensional representation of 
a closed subgroup $H$ of $G$ and $C^\infty(G;V_\varpi)$ denote the space of 
$V_\varpi$-valued $C^\infty$-functions on $G$.
Then the space of $C^\infty$-sections $C^\infty(G/H;\varpi)$ of 
the $G$-homogeneous bundle associated
to $\varpi$ is
\[
 \{f\in C^\infty(G;V_\varpi)\,;\,f(gh)=\varpi^{-1}(h)f(g)
  \text{ for }\forall h\in H\}.
\]

Consider the case when $H=K$.
Because of the decomposition $G=KAK$ the function 
$f\in C^\infty(G/K;\varpi)$ is determined by its restriction on $KA$ and by
the natural map $K\times A\to KA$ the restriction can be considered as
a function $\bar f$ on $K\times A$.  Then the action of the differential 
operator $L_2$ to $\bar f$ is
\begin{align*}
 \bar L_2 &=   \sum_{1\le i<j\le n}
       \Bigl((\tilde\vartheta_{i-1}-\tilde\vartheta_i)
       (\tilde\vartheta_{j-1}-\tilde\vartheta_j)
       -\frac{t_{ij}^2}{1-t_{ij}^2}(\vartheta_{i-1}-\vartheta_i
        - \vartheta_{j-1} + \vartheta_j)\\
  &\quad
   -\frac{t_{ij}^2}{(1-t_{ij}^2)^2}
            \bigl(\Ad(a^{-1})(E_{ij}-E_{ji})^2
            +\varpi(E_{ij}-E_{ji})^2\bigr)\notag\\
   &\quad
  - \frac{t_{ij}(1+t_{ij}^2)}{(1-t_{ij}^2)^2}\Ad(a^{-1})(E_{ij}-E_{ji})\cdot
   \varpi(E_{ij}-E_{ji})\Bigr)
\end{align*}
at $(k,a) \in K\times A$, which follows from \eqref{eq:Ktype}.
Here the induced representation of the Lie algebra $\mathfrak k$ of $K$ 
is also denoted by $\varpi$.

Let $(\delta,V_\delta)$ be an irreducible representation of $K$.
Then the $\delta$-component of $C^\infty(G/K;\varpi)$ is an element
$f\in V\otimes C^\infty(G/K;\varpi)$ which satisfies 
\[
 \frac{d}{dt}f(e^{tX}g)\Bigm|_{t=0}=\bigl(\delta(X)f\bigr)(g)
\]
for $X\in\mathfrak k$.
Hence the function $f$ is determined by its restriction $\bar f$ on $A$ and 
the action of the operator $L_2$ to $\bar f$ is
\begin{align}
  \bar L_2 &=   \sum_{1\le i<j\le n}
       \Bigl((\tilde\vartheta_{i-1}-\tilde\vartheta_i)
       (\tilde\vartheta_{j-1}-\tilde\vartheta_j)
       -\frac{t_{ij}^2}{1-t_{ij}^2}(\vartheta_{i-1}-\vartheta_i
       -\vartheta_{j-1} + \vartheta_j)\notag\\
  &\quad
   -\frac{t_{ij}^2}{(1-t_{ij}^2)^2}
            \bigl(\delta(E_{ij}-E_{ji})^2+\varpi(E_{ij}-E_{ji})^2\bigr)
   \label{eq:Sph}\\
   &\quad
  -\frac{t_{ij}(1+t_{ij}^2)}{(1-t_{ij}^2)^2}\delta(E_{ij}-E_{ji})
  \otimes\varpi(E_{ij}-E_{ji})
  \Bigr).\notag
\end{align}
Note that the operator $P=\bar L_2$ satisfies the assumption of 
Corollary~\ref{cor:commute}.

When $G$ is $SL(2,\mathbb R)$ or its universal covering group and $\bar f$
is an eigenfunction of $L_2$, we can put $\varpi(E_{12}-E_{21})=\sqrt{-1}k$ and $\delta(E_{12}-E_{21})=-\sqrt{-1}m$ for certain numbers $k$ and $m$ and
\[
  \Bigl(\vartheta^2+\frac14-\frac{1+t^2}{1-t^2}\vartheta
  +\frac{t(k-mt)(m-kt)}{(1-t^2)^2}-(\lambda+\frac12)^2
  \Bigl)\bar f=0.
\]
Put $t=e^{-x}$ and $u=\bar f$. 
Then $\vartheta=-\frac{d}{dx}$ and
\begin{gather*}
 u''+\coth x\cdot u'-\frac{(k+m)^2}{4\sinh^2x}u+
 \frac{km}{4\sinh^2\frac x2}u=\lambda(\lambda+1)u,\\
  \frac{d^2 v}{d z^2} - \frac{(k+m-1)(k+m+1)}{\sinh^22z} v 
     + \frac{km}{\sinh^2z}v = \bigl(2\lambda+1\bigr)^2 v.
\end{gather*}
by denoting $v=\sinh^{\frac12}x\cdot u$ and $z=\frac x2$.

Then for $\tilde v=\sinh^{m}z\cdot \sinh^{-\frac{k+m+1}2}2z\cdot v$ and
$w=-\sinh^2 z$ we have
\[
 \begin{split}
  w(1-w)\frac{d^2\tilde v}{d w^2}
  + \Bigl(\frac{k-m+2}2 - (k+2)w\Bigr)\frac{d\tilde v}{dw}
  -\Bigl(\frac k2-\lambda\Bigr)\Bigl(\frac k2+\lambda + 1\Bigr)\tilde v=0
 \end{split}
\]
and hence $\bar f$ is a linear combination of the functions
\[
 \begin{cases}
 \sinh^{\frac{k-m}2}z\cdot\cosh^{\frac{k+m}2}z\cdot 
F(\frac k2- \lambda,\frac k2+\lambda+1,
\frac{k-m}2+1;-\sinh^2 z),\\
 \sinh^{\frac{m-k}2}z\cdot\cosh^{\frac{k+m}2}z\cdot F(
 \frac m2-\lambda,\frac m 2+\lambda+1,
 \frac{m-k}2+1;-\sinh^2 z\bigr).
 \end{cases}
\]
Thus it is clear that the non-zero real analytic solution $\bar f$ defined 
in a neighborhood of the point $z=0$ exists if and only if $k-m\in2\mathbb Z$.
Here $F(\alpha,\beta,\gamma;z)$ denotes the Gauss hypergeometric function
(cf.~\cite{WW}).

Next we assume $H=N$ and $\varpi$ is a character of $N$.
Then there exist complex numbers $c_1,\ldots,c_{n-1}$ such that
\[
 \varpi(e^{\sum_{1\le i<j\le n}s_{ij}E_{ij}})
 = r^{\sqrt{-1}(c_1s_{12}+\cdots+c_{n-1}s_{n-1,n})}.
\]
The element $f\in C^\infty(G/N;\varpi)$ is determined by the restriction 
$\bar f=f|_{KA}$ and it follows from \eqref{eq:Ntype} that 
the operation of $L_2$ to $\bar f$ is 
\begin{equation}\label{eq:KWhit}
 \sum_{1\le i<j\le n}
  (\tilde\vartheta_{i-1}-\tilde\vartheta_i)
  (\tilde\vartheta_{j-1}-\tilde\vartheta_j)
  +\sum_{1\le i<n}
  \bigl(c_i^2t_i^2 +\sqrt{-1} c_it_i(E_{i,i+1}-E_{i+1,i})\bigr).
\end{equation}
Hence if $G=SL(2,\mathbb R)$, the eigenfunction $f$ of $L_2$ of the 
$\delta$-component of 
$C^\infty(G/N;\varpi)$ with $\delta(E_{12}-E_{21})=\sqrt{-1}m$ satisfies 
\[
  \Bigl(-(\vartheta-\tfrac12)^2+c_1^2t^2-c_1mt+(\lambda+\tfrac12)^2\Bigr)f|_A = 0
\]
and hence
\[
  \frac{d^2}{dt^2}(f|_A) - \Bigl(c_1^2 
  - \frac{c_1m}t+\frac{\lambda(\lambda+1)}{t^2}\Bigr)(f|_A)=0.
\]
If we put $u(x)= e^{\frac x2}\bigl(f|_A(e^{-x})\bigr)$, then
\[
  u'' - \bigl(c_1e^{-2x} - c_1me^{-x}\bigr)u = \lambda(\lambda+1)u.
\]
Denoting $W(\pm2c_1t)=w(t)$, we have the Whittaker equation (cf.~\cite{WW})
\[
  W'' + \Bigl(-\frac14\pm\frac{m}{2t}
  +\frac{\frac14-(\lambda+\frac12)^2}{t^2}\Bigr)W=0.
\]
\section{Completely integrable quantum systems}\label{sec:CI}
A Schr\"odinger operator
\begin{equation*}
   P = \displaystyle\sum_{k=1}^n\frac{\p^2}{\p x_k^2} + R(x_1,\ldots,x_n)
\end{equation*}
of $n$ variables is called 
\emph{completely integrable} if there exist $n$ algebraically independent
differential operators $P_k$ such that
\begin{equation*}
  [P_i,P_j]= 0\quad\text{for \ }1\le i<j\le n
  \quad\text{and \ }P\in\mathbb C[P_1,\ldots,P_n].
\end{equation*}
Under the coordinate system $(t_1,\ldots,t_n)$ with
\begin{equation*}
   t_1=e^{x_1-x_2},\ldots,t_{n-1}=e^{x_{n-1}-x_n}, t_n=e^{x_n},
\end{equation*}
the Schr\"odinger operators $P$
which belong to $\mathcal D_*$ and have elements $Q\in\mathcal D_*$ satisfying
\begin{equation*}
  Q = \sum_{k=1}^n \frac{\p^4}{\p x_k^4} + Q'\quad\text{with \ }\ord Q'<4
\end{equation*}
are classified in \cite{OI2} and proved to be completely integrable 
(cf.\ \cite{OI1} and \cite{OI3}).
They are reduced to the Schr\"odinger operators
with the potential functions $R(x_1,\ldots,x_n)$ in the following list.
\begin{gather*}
\allowdisplaybreaks
 \begin{aligned}
  &
  \sum_{1\le i<j\le n}C_1\bigl(\sh^{-2}\tfrac{x_i+x_j}2
  +\sh^{-2}\tfrac{x_i-x_j}2\bigr)\\
  &\quad+\sum_{k=1}^n\bigl(
    C_2\sh^{-2}x_k + C_3\sh^{-2}\tfrac{x_k}2\bigr),
  \end{aligned}
\tag{\rm Trig-$BC_n$-reg}
\allowdisplaybreaks\\
 \begin{aligned}
 &\sum_{1\le i<j\le n}C_1\sh^{-2}\tfrac{x_i-x_j}2
 + \sum_{k=1}^n\bigl(C_2e^{x_k} + 
   C_3e^{2x_k}\bigr),
 \end{aligned}
\tag{\rm Trig-$A_{n-1}$-bry-reg}
\allowdisplaybreaks\\
 \begin{aligned}
 &C_1\sum_{i=1}^{n-1}e^{x_i-x_{i+1}}
 + C_1e^{x_{n-1}+x_n}
 + C_2\sh^{-2}\tfrac{x_n}2 +  C_3\sh^{-2}x_n,
 \end{aligned}
\tag{\rm Toda-$D_n$-bry}
\allowdisplaybreaks\\
 \begin{aligned}
 &C_1\sum_{i=1}^{n-1}e^{x_i-x_{i+1}}
 + C_2e^{x_n}+ C_3e^{2x_n}.
 \end{aligned}
\tag{\rm Toda-$BC_n$}
\end{gather*}
Here $C_1$, $C_2$ and $C_3$ are any complex numbers.

We can generalize the Schr\"odinger operators 
in terms of root systems (cf.~\cite{OP}).
Let $\Sigma$ be an irreducible root system with rank $n$,
$\Sigma^+$ a positive system of $\Sigma$
and $\Psi\subset\Sigma$ a fundamental system of $\Sigma^+$.
Then
$\Sigma$ is identified with a finite subset of a Euclidean space $\mathbb R^n$
and
\begin{align}\label{eq:roottrig}
 P = \sum_{k=1}^n\frac{\p^2}{\p x_k^2}
     + \sum_{\alpha\in\Sigma^+}
       \frac{C_\alpha}{\sh^2\frac{\langle\alpha,x\rangle}2}
\quad(C_\alpha\in\mathbb C,\ C_\alpha=C_\beta\text{ if }|\alpha|=|\beta|)
\end{align}
and
\begin{align}\label{eq:rootToda}
 P = \sum_{k=1}^n\frac{\p^2}{\p x_k^2}
     + \sum_{\alpha\in\Psi}e^{\langle\alpha,x\rangle}
\end{align}
are Schr\"odinger operators of Heckman-Opdam's hypergeometric system
(cf.\ \cite{HO}) and Toda finite chain (cf.\ \cite{To}) corresponding 
to the fundamental system $\Psi$, respectively.
They are in $\mathcal D_*$ under the coordinate system
\begin{equation*}
  t_k = e^{\langle\alpha_k,x\rangle}\quad\text{for \ }k=1,\ldots,n
\end{equation*}
with $\Psi=\{\alpha_1,\ldots,\alpha_n\}$ and known to be completely integrable.

If $\Sigma$ is of type $BC_n$, then
\begin{equation*}
 \Sigma^+=\{e_i-e_j,\,e_k,\,2e_k\,;\,1\le i<j\le n,\, 1\le k\le n \}
\end{equation*}
and the Schr\"odinger operators \eqref{eq:roottrig} and
\eqref{eq:rootToda} correspond to (Trig-$BC_n$-reg) or (Toda-$BC_n$).
If $\Sigma$ is of other classical type, the operators also correspond to 
special cases of (Trig-$A_{n-1}$-bry) or (Toda-$D_n$-bry) or (Toda-$BC_n$).

The potential functions $R(x)$ of known completely integral quantum systems
which may not have regular singularities at infinity 
are expressed by functions with one variable.
If $P_2$ and $P_3$ are operators of order 4 and 6 with 
the highest order terms
$\sum_{k=1}^n\frac{\p^4}{\p x_k^4}$ and $\sum_{k=1}^n\frac{\p^6}{\p x_k^6}$,
respectively, this is proved by \cite{Wa} in general.
We will examine this in the case when $n=2$.

\begin{theorem}\label{thm:CI2}
Let $\ell$ be a positive integer.
Suppose the differential operators
\begin{equation}\label{eq:com2v}
 \begin{split}
  P &= \frac{\p^2}{\p x^2}+\frac{\p^2}{\p y^2}+R(x,y),\\
  Q &= \sum_{i=0}^m c_i\frac{\p^m}{\p x^{m-i}\p y^i}
       + \sum_{i+j\le m-2} S_{i,j}(x,y)\frac{\p^{i+j}}{\p x^i\p y^j}
\end{split}
\end{equation}
satisfy $[P,Q]=0$ and $\sigma_m(Q)\notin\mathbb C[\sigma(P)]$.
Here $R(x,y)$ and $S_{i,j}(x,y)$ are square matrices of size $\ell$ whose
components are functions of $(x,y)$ and $c_i\in\mathbb C$. 
Put
\begin{equation}\label{eq:com2v2}
 \Bigl(\xi\frac{\p}{\p\tau}-\tau\frac{\p}{\p\xi}\Bigr)
  \sum_{i=0}^mc_i\xi^{m-i}\tau^i
 =\prod_{\nu=1}^L(a_\nu\xi - b_\nu\tau)^{m_\nu}
\end{equation}
with suitable $(a_\nu,b_\nu)\in\mathbb C^2\setminus\{0\}$ satisfying
${a_\nu}{b_\mu}\ne {a_\mu}{b_\nu}$ for $\mu\ne\nu$.
Here $m_\nu$ are positive integers and $m_1+\cdots+m_L=m$.
Then
\begin{equation}
 R(x,y) = \sum_{\nu=1}^L\sum_{i=0}^{m_\nu-1}
 (b_\nu x + a_\nu y)^i R_{\nu,i}(a_\nu x-b_\nu y)
\end{equation}
with $m$ square matrices of size $\ell$ whose components are functions 
$R_{\nu,i}(t)$ of the one variable $t$.
\end{theorem}
\begin{proof}
The coefficients of $\frac{\p^{m+1}}{\p x^{m-1-j}\p y^j}$ 
in the expression $[P,Q]$ for \eqref{eq:com2v} show
\begin{equation*}
 2\p_xS_{m-2-j,j}+2\p_yS_{m-1-j,j-1}=c_j(m-j)\p_xR + c_{j+1}(j+1)\p_yR
\end{equation*}
for $j=0,\ldots,m-1$. Hence the theorem follows from the following equation.
\begin{align*}
  0&=2\sum_{j=0}^{m-1}(-1)^j\bigl(\p_x^j\p_y^{m-j} S_{m-2-j,j}
  +\p_x^{j+1}\p_y^{m-1-j} S_{m-1-j,j-1}\bigr)\\
   &=2\sum_{j=0}^{m-1}(-1)^j\p_x^j\p_y^{m-1-j}
    \bigl(c_j(m-j)\p_xR + c_{j+1}(j+1)\p_yR\bigr)\\
   &=\sum_{j=0}^{m-1}(-1)^jc_j(m-j)\p_x^{j+1}\p_y^{m-1-j}R
     +\sum_{j=0}^{m-1}(-1)^jc_{j+1}(j+1)\p_x^{j}\p_y^{m-j}R\\
   &=\Bigl(\bigl(\xi\frac{\p}{\p\tau}-\tau\frac{\p}{\p\xi}\bigr)\sum_{i=0}^mc_i\xi^{m-i}\tau^i\Bigr)\Bigm|_{\xi=\p_y,\,\tau=-\p_x}R.
\end{align*}
\end{proof}


\begin{thebibliography}{CFV}
\bibitem[Ca]{Ca} Capelli, A.:
\"Ueber die Zur\"uckf\"uhrung der Cayley'schen Operation
$\Omega$ auf gew\"ohnliche Polar-Operationen.
Math.\ Ann.\ {\bf 29}, 331-338(1887).

\bibitem[Ha]{Ha} Harish-Chandra:
Collected Papers I-IV. Springer(1989).

\bibitem[HO]{HO} Heckman, G.J., Opdam, E.M.:
Root system and hypergeometric functions. I.
Comp.\ Math.\ {\bf 64}, 329--352(1987).

\bibitem[K--]{K--} Kashiwara, M., Kowata, A., Minemura, K.,  Okamoto, K,, 
Oshima, T.,  Tanaka, M.:
Eigenfunctions of invariant differential operators on a symmetric space.
Ann.\ of Math.\ {\bf 107}, 1--39(1978).

\bibitem[KO]{KO} Kashiwara, M., Oshima, T.:
Systems of differential equations with regular singularities and their 
boundary value problems.
Ann.\ of Math.\ {\bf 106}, 145--200(1977).

\bibitem[OP]{OP} Olshanetsky, M.A.,  Perelomov, A.M.:,
Quantum integrable systems related to Lie algebras.
Phys.\ Rep.\ {\bf 94}, 313--404(1983).

\bibitem[O1]{OV0} Oshima, T.:
A boundary value problem on a Riemannian symmetric space.
S\-uriKaisekiKenky\-usho Kouky\-uroku {\bf 249}, 10--21(1975),
in Japanese.

\bibitem[O2]{OR1} Oshima, T.:
A realization of Riemannian symmetric spaces.
J.\ Math.\ Soc.\ Japan {\bf 53}, 117--132(1978).

\bibitem[O3]{OV1} Oshima, T.:
A definition of boundary values of solutions of partial differential equations 
with regular singularities.
Publ.\ RIMS Kyoto Univ.\ {\bf 19}, 1203--1230(1983).

\bibitem[O4]{OV2} Oshima, T.: 
Boundary value problems for systems of linear partial differential equations 
with regular singularities. 
Advanced Studies in Pure Math.\ {\bf 4}, 391--432(1984).

\bibitem[O5]{O-Asymp} Oshima, T.:
Asymptotic behavior of spherical functions on semisimple symmetric spaces.
Advanced Studies in Pure Math.
\textbf{14}, 561--601 (1988).

\bibitem[O6]{OR2} Oshima, T.: 
A realization of semisimple symmetric spaces and construction of boundary 
value maps.
Advanced Studies in Pure Math.\ {\bf 14}, 603--650(1988).

\bibitem[O7]{OI1} Oshima, T.: 
Completely integrable systems with a symmetry in coordinates.
Asian Math.\ J.\ {\bf 2}, 935--956(1998).

\bibitem[O8]{OI2} Oshima, T.: 
A class of completely integrable quantum systems associated with 
classical root systems.
Indag.\ Mathem.\ {\bf 16}, 655--677(2005).

\bibitem[O9]{OI3} Oshima, T.: 
Completely integrable quantum systems associated 
with classical root systems.
preprint, 2005, 45pp, math-ph/0502028.

\bibitem[OS]{OS} Oshima, T., Sekiguchi, J.:
Eigenspaces of invariant differential operators on an affine symmetric space.
Invent.\ Math.\ {\bf 57}, 1--81(1980).

\bibitem[Su]{Su} Sutherland, B:
Exact results for a quantum many-body problem in one dimension II.
Phys.\ Rev.\ {\bf A5}, 1372--1376(1972).

\bibitem[To]{To} Toda, M.:
Wave propagation in anharmonic lattice.
J.\ Phy.\ Soc.\ Japan {\bf 23}, 501--506(1967).

\bibitem[Wa]{Wa} Wakida, S.:
Quantum integrable systems associated with classical Weyl groups.
MA Thesis, University of Tokyo, Tokyo(2004).

\bibitem[W]{WW} Whittaker, E.T., Watson, G.N.:
A Course of Modern Analysis, Fourth Edition. 
Cambridge University Press (1927).
\end{thebibliography}
\end{document}